\documentclass[leqno,12pt]{article}
\usepackage{amsmath,amssymb,color,a4wide}
\usepackage{MnSymbol}
\usepackage{tikz-cd}
\usepackage{xr-hyper}
\usepackage{hyperref}

\usepackage[curve,matrix,arrow,cmtip]{xy}
\usepackage{verbatim}
\NoComputerModernTips

\newdir^{ (}{{}*!/-3pt/\dir^{(}}
\newdir_{ (}{{}*!/-3pt/\dir_{(}}

\def\BAfF{{\BA^{\rm f}_F}}
\def\bigtimes{\mathop{\raisebox{-2pt}{\huge $\times$\kern-1pt}}}
\def\BGa{{\BG_{\mathrm{a}}}}
\def\BGacomma#1{{\BG_{\mathrm{a},#1}}}
\def\hatA{\smash{\hat{A}}}

\def\OM{\mathchoice
  {\rlap{\kern3.2pt$\overline{\phantom{L}}$}M}
  {\rlap{\kern3.2pt$\overline{\phantom{L}}$}M}
  {\rlap{\kern2.4pt$\scriptstyle\overline{\phantom{L}}$}M}
  {\rlap{\kern1.8pt$\scriptscriptstyle\overline{\phantom{L}}$}M}}

\let\le\leqslant
\let\ge\geqslant

\def\an{{\rm an}}
\def\alg{{\rm alg}}

 \def\End{\mathop{\rm End}\nolimits}

 \def\Spec{\mathop{\rm Spec}\nolimits}
 
 \def\deg{\mathop{\rm deg}\nolimits}

\def\Lie{\mathop{\rm Lie}\nolimits}

\def\GL{\mathop{\rm GL}\nolimits}
\def\SL{\mathop{\rm SL}\nolimits}

\def\id{{\rm id}}
\def\Id{{\rm Id}}

\let\phi\varphi
\let\theta\vartheta
\let\epsilon\varepsilon
\let\setminus\smallsetminus

\newtheorem{Thm}{Theorem}[section]
\newtheorem{Prop}[Thm]{Proposition}
\newtheorem{Lem}[Thm]{Lemma}

\newtheorem{Def}[Thm]{Definition}

\newtheorem{Rem}[Thm]{Remark}

\def\UseTheoremCounterForNextEquation{\setcounter{equation}{\value{Thm}}\addtocounter{Thm}{1}}

\def\qed{{\hskip0pt\unskip\unskip\nobreak\hfil\penalty50
          \hskip1em\hbox{}\nobreak\hfil
           {$\square$}
          \parfillskip=0pt\finalhyphendemerits=0
          \par}\medskip}

\newenvironment{Proof}
               {\noindent{\bf Proof.}\ }
               {\qed}
							
               {\noindent{\bf Proof Sketch.}\ }
               {\qed}

               {\noindent{\bf Proof of #1.}\ }
               {\qed}

\newcommand{\BA}{{\mathbb{A}}}

\newcommand{\BC}{{\mathbb{C}}}

\newcommand{\BF}{{\mathbb{F}}}
\newcommand{\BG}{{\mathbb{G}}}

\newcommand{\Fp}{{\mathfrak{p}}}

\newcommand{\CL}{{\cal L}}
\newcommand{\CM}{{\cal M}}

\newcommand{\CO}{{\cal O}}

\newcommand{\CS}{{\cal S}}
\newcommand{\CT}{{\cal T}}
\newcommand{\CU}{{\cal U}}
\newcommand{\CV}{{\cal V}}
\newcommand{\CW}{{\cal W}}

\newbox\mybox
\def\arrover#1{\mathrel{
       \setbox\mybox=\hbox spread 1.4em
              {\hfil$\scriptstyle#1$\hfil}
       \vbox{\offinterlineskip\copy\mybox
             \hbox to\wd\mybox{\rightarrowfill}}}}

\def\larrover#1{\mathrel{
       \setbox\mybox=\hbox spread 1.4em
              {\hfil$\scriptstyle#1\vphantom{g}$\hfil}
       \vbox{\offinterlineskip\copy\mybox
             \hbox to\wd\mybox{\leftarrowfill}}}}

\def\ontoover#1{\mathrel{
       \setbox\mybox=\hbox spread 1.4em
              {\hfil$\scriptstyle#1\vphantom{g}$\hfil}
       \vbox{\offinterlineskip\copy\mybox
             \hbox to\wd\mybox{\rightarrowfill\hskip-2.8mm
                               $\rightarrow$}}}}
\def\leftontoover#1{\mathrel{
       \setbox\mybox=\hbox spread 1.4em
              {\hfil$\scriptstyle#1\vphantom{g}$\hfil}
       \vbox{\offinterlineskip\copy\mybox
             \hbox to\wd\mybox{$\leftarrow$\hskip-2.8mm
                               \leftarrowfill}}}}
\let\longto\longrightarrow
\let\into\hookrightarrow
\let\onto\twoheadrightarrow

\def\Cinf{{\BC}_\infty}
\def\Finf{F_\infty}

\def\Bigskip{\bigskip\bigskip}

\newbox\dotDdbox
\newbox\dotDtbox
\newbox\dotDsbox
\newbox\dotDssbox
\setbox\dotDdbox=\vbox{\hbox{\kern3.5pt\bf.}\vskip-4.5pt\hbox{$D$}}
\setbox\dotDtbox=\vbox{\hbox{\kern3.5pt\bf.}\vskip-4.5pt\hbox{$D$}}
\setbox\dotDsbox=\vbox{\hbox{\kern2.1pt.}
                       \vskip-6.6pt\hbox{$\scriptstyle D$}}
\setbox\dotDssbox=\vbox{\hbox{\kern1.6pt.}
                        \vskip-8.1pt\hbox{$\scriptscriptstyle D$}}

\newcommand{\Malg}{\CM^\alg}

\newcommand{\PropositionCite}{Prop.\,}
\newcommand{\PropositionsCite}{Props.\,}
\newcommand{\TheoremCite}{Thm.\,}
\newcommand{\DefinitionCite}{Def.\,}

\begin{document}

\title{Drinfeld modular forms of arbitrary rank\\
Part II: Comparison with Algebraic Theory}

\author{Dirk Basson \and Florian Breuer$^{1,2}$ \and Richard Pink$^2$}

\footnotetext[1]{Supported by the Alexander von Humboldt foundation, and by the NRF grant BS2008100900027.}
\footnotetext[2]{Supported through the program ``Research in Pairs'' by Mathematisches Forschungsinstitut Oberwolfach in 2010.}

\date{May 27, 2018}
\maketitle

\centerline{To Oliver}
\bigskip

\begin{abstract}
This is the second of a series of articles providing a foundation for the theory of Drinfeld modular forms of arbitrary rank. In the present part, we compare the analytic theory with the algebraic one that was begun in a paper of the third author. For any arithmetic congruence subgroup and any integral weight we establish an isomorphism between the space of analytic modular forms with the space of algebraic modular forms defined in terms of the Satake compactification. From this we deduce the important result that this space is finite dimensional.
\end{abstract}

{\advance\baselineskip by -6pt
\tableofcontents
}


\externaldocument{Part_1}
\externaldocument{Part_3}
\setcounter{section}{6}

\Bigskip


\newpage
\noindent {\Large\bf Introduction}
\bigskip

This is part II of a series of articles together with \cite{BBP1} and \cite{BBP3}, whose aim is to provide a foundation for the theory of Drinfeld modular forms of arbitrary rank. Part I developed the basic analytic theory, including $u$-expansions and holomorphy at infinity. In the present Part II we identify the analytic modular forms from Part I with the algebraic modular forms defined in \cite{Pink} and deduce qualitative consequences such as the finite dimensionality of the space of modular forms of given level and weight. Part III will illustrate the general theory by constructing and studying some important families of modular forms.

\medskip
By definition, \emph{weak Drinfeld modular forms} of weight $k$ are holomorphic functions on the rigid analytic Drinfeld period domain $\Omega^r$ that satisfy a certain twisted transformation law under the action of an arithmetic congruence subgroup $\Gamma<\GL_r(F)$. 
\emph{Drinfeld modular forms} are weak Drinfeld modular forms that are holomorphic at infinity after transformation by all elements of $\GL_r(F)$. By construction these seem to be purely analytic objects, but in this article we identify them with objects from algebraic geometry, as follows.

Roughly speaking, the quotient $\Gamma\backslash\Omega^r$ is the set of $\Cinf$-valued points of a certain moduli space of Drinfeld modules~$M$, which is an algebraic variety over~$\Cinf$. The transformation law means that weak modular forms of weight $k$ can be interpreted as holomorphic sections of $\CL^k$ for a certain invertible sheaf $\CL$ on~$M$, at least if $\Gamma$ is sufficiently small. Here $\CL$ is the dual of the relative Lie algebra of the universal Drinfeld module over~$M$. Since $M$ is affine of dimension~$r-1$, for $r\ge2$ there is an abundance of non-algebraic holomorphic sections of~$\CL^k$. (So the analogue of the K\"ocher principle for Siegel modular forms does not hold.) 

To algebraise Drinfeld modular forms, we translate the condition at infinity into a condition on a compactification $\OM$ of the moduli space~$M$. For this we use the Satake compactification that was constructed analytically by Kapranov \cite{Kapranov} in the special case $A=\BF_q[t]$ and by H\"aberli \cite{Haeberli} in general, and algebraically by the third author in \cite{Pink}. By \cite{Pink} the sheaf $\CL$ extends naturally to an invertible sheaf on~$\OM$, again denoted~$\CL$, which is constructed as the dual of the relative Lie algebra of the unique generalised Drinfeld module over~$\OM$ that extends the universal Drinfeld module over~$M$. 

The main result of this article, Theorem \ref{AAMFThm1}, states that the analytic Drinfeld modular forms of weight $k$ correspond precisely to the sections of $\CL^k$ over~$\OM$. Since $\OM$ is a projective algebraic variety, it follows that the space of modular forms of each weight $k$ is finite dimensional, and that the graded ring of modular forms of all weights for fixed~$\Gamma$ is a normal integral domain that is finitely generated as a $\Cinf$-algebra: see Theorem \ref{FiniteDimension}.
In the case $r=2$ all this was done in Goss's thesis \cite{GossES}.

\medskip
Establishing these results with adequate precision requires a fair amount of technical details. 
For later use we also discuss the action of $\GL_r(F)$ as well as Hecke operators. 
As this article belongs to a whole series with \cite{BBP1} and \cite{BBP3}, we  number the sections of all parts consecutively. Thus Sections 1--6 appear in Part I and Sections 13--17 in Part III. All the definitions and notation from Part I remain in force, and we refer to proclamations in the other parts without any special indication.

\subsubsection*{Outline  of this paper}

As a preparation for the modular interpretation of $\Gamma\backslash\Omega^r$, in Section \ref{Sec:UFDM} we construct the universal family of Drinfeld modules over $\Omega^r$ and its level structures. We also study its behaviour at the standard boundary component. In Proposition \ref{GenDrinModCodim1} we show that the universal family descends to a family over $\Gamma_U\backslash\Omega^r$ which extends naturally to a generalised Drinfeld module over the larger domain $\CU$ obtained by adjoining a copy of~$\Omega^{r-1}$.

In Section \ref{Sec:DMS} we construct the precise identification of $\Gamma\backslash\Omega^r$ with 
a moduli space of Drinfeld modules. This requires working with the ring of finite ad\`eles $\BAfF$ of~$F$ and identifying $\Gamma\backslash\Omega^r$ with a connected component of a double quotient of the form 
$$\GL_r(F)\!\bigm\backslash\! \bigl( \Omega^r\times\GL_r(\BAfF)/K \bigr)$$
for an open compact subgroup $K<\GL_r(\hatA)$. That in turn can be identified naturally with the space of $\Cinf$-valued points $M^r_{A,K}(\Cinf)$ on a certain algebraic moduli space of Drinfeld modules $M^r_{A,K}$. This identification requires a precise description of the universal family and its level structure. Working ad\`elically also entails that $M^r_{A,K}$ is an algebraic variety over the given global field $F$ itself, which eventually shows that the space of modular forms for $\Gamma$ comes from a vector space over a certain finite abelian extension of $F$ instead of~$\Cinf$. 

As explained in Remark \ref{ConventionsRemA}, there are different conventions about whether $\Omega^r$ consists of row or column vectors and about how $\GL_r(\Finf)$ acts on it. In this series of articles we have chosen to use column vectors and left multiplication. This affects the way that the universal family of Drinfeld modules on $\GL_r(F)\backslash(\Omega^r\times\GL_r(\BAfF)/K)$ must be described. As our convention differs from that of \cite{Pink}, several formulas from there have to be transformed to be used here. For instance, in the isomorphism (\ref{AnalUniv1}) a double coset $[(\omega,g)]$ now corresponds to a point on the moduli space that was represented by the double coset $[(\omega^T,(g^T)^{-1})]$ in \cite{Pink}. The change in convention also affects the functoriality in Proposition \ref{Func1}, in whose proof the precise relationship is indicated.
We wish to apologise for the resulting inconvenience.

In Section \ref{Sec:Satake} we review the relevant facts about the Satake compactification of $\OM^r_{A,K}$ of $M^r_{A,K}$. The crucial properties in Proposition \ref{Codim1Prop1} are that the composite map $\Gamma_U\backslash\Omega^r \onto \Gamma\backslash\Omega^r \into M^r_{A,K}(\Cinf)$ extends to an \'etale morphism $\CU\to\OM^r_{A,K}(\Cinf)$ for the larger domain $\CU$ from Section \ref{Sec:UFDM} and that, repeating this after transformation by arbitrary elements of $\GL_r(\BAfF)$, the images of these maps cover a Zariski open subset $M^{r,+}_{A,K}(\Cinf)$ of $\OM^r_{A,K}(\Cinf)$ whose closed complement has codimension $\ge2$. Using this map we can identify the pullback of the generalised Drinfeld module on $\OM^r_{A,K}$ with that constructed over $\CU$ in Section \ref{Sec:UFDM}.

In Section \ref{Sec:AAMF} we use these facts to show that an analytic modular form is holomorphic at infinity if and only if the corresponding section of $\CL^k$ over $M^r_{A,K}(\Cinf)$ extends holomorphically to a section over $M^{r,+}_{A,K}(\Cinf)$. By rigid analytic analogues of the Hartogs principle and of GAGA the latter condition is equivalent to being the restriction of a section of $\CL^k$ over $\OM^r_{A,K}(\Cinf)$ in the algebro-geometric sense, thereby establishing our first main result, Theorem \ref{AAMFThm1}.

This earns us our piece of cake in Section \ref{Sec:Fin}, where we deduce that the space of modular forms of each weight $k$ is finite dimensional, and that the graded ring of modular forms of all weights for fixed~$\Gamma$ is a normal integral domain that is finitely generated as a $\Cinf$-algebra.

The final Section \ref{HeckeOps} explains how the comparison isomorphism between analytic and algebraic modular forms behaves under Hecke operators on both sides.



\section{Universal family of Drinfeld modules}
\label{Sec:UFDM}

As a preparation for the following sections, we construct the universal family of Drinfeld modules on 
$\Gamma\backslash\Omega^r$ associated to an $A$-lattice $L\subset F^r$ and study its behaviour at the standard boundary component. We first review the necessary details about Drinfeld modules and generalised Drinfeld modules. 

\medskip
Consider any scheme $S$ over~$F$. For any line bundle $E$ on~$S$, let $\End_{\BF_q}(E)$ denote the ring of $\BF_q$-linear endomorphisms of the group scheme underlying~$E$. 
(These endomorphisms need not commute with scalar multiplication by $\CO_S$.) 
By \cite[\S5]{Drinfeld1}, any such endomorphism is a finite sum $\sum_{i}b_i\tau^i$ for sections $b_i\in H^0(S,E^{1-q^i})$, where $\tau : {E\to E^q}$, ${x\mapsto x^q}$ denotes the $q$-power Frobenius morphism.
Set $\deg(a) := \dim_{\BF_q}(A/(a))$ for any $a\in A\setminus\{0\}$ and $\deg(0):=-\infty$.

Recall that a \emph{Drinfeld $A$-module of rank $r$ over~$S$} is a pair $(E,\varphi)$ consisting of a line bundle $E$ over $S$ and a ring homomorphism
\UseTheoremCounterForNextEquation
\begin{equation}\label{DrinfeldModuleDef}
\varphi:\ A\to\End_{\BF_q}(E),\ \ 
a\mapsto \varphi_a=\!\!\sum_{i=0}^{r\deg(a)}\!\! \varphi_{a,i}\tau^i
\end{equation}
with $\varphi_{a,i}\in H^0(S,E^{1-q^i})$ satisfying the two conditions:
\begin{enumerate}
\item[(a)] The derivative $d\varphi:a\mapsto \varphi_{a,0}$ is the structure homomorphism $A\into F\to H^0(S,\CO_S)$.
\item[(b)] For any $a\in A\setminus\{0\}$ the term $\varphi_{a,r\deg(a)}$ is a nowhere vanishing section of $E^{1-q^{r\deg(a)}}$.
\end{enumerate}
If instead of (b) we require only:
\begin{enumerate}
\item[(c)] For any point $s\in S$ and any non-constant $a\in A$ there exists $i>0$ with $\varphi_{a,i}\not=0$;
\end{enumerate}
we obtain the notion of a \emph{generalised Drinfeld $A$-module of rank $\le r$ over $S$} from \cite[\DefinitionCite 3.1]{Pink}.
Over any point $s\in S$, the map $\varphi$ then defines a Drinfeld $A$-module of some rank $r_s$ satisfying $1\le r_s\le r$.

An \emph{isomorphism} of (generalised or not) Drinfeld $A$-modules over $S$ is an isomorphism of line bundles that is equivariant with respect to the action of $A$ on both sides.
Furthermore, following \cite[\DefinitionCite 3.8]{Pink}, a generalised Drinfeld $A$-module $(E,\phi)$ over $S$ is called \emph{weakly separating} if, for any Drinfeld $A$-module $(E',\phi')$ over any field $F'$ containing $F$, at most finitely many fibers of $(E,\phi)$ over $F'$-valued points of $S$ are isomorphic to $(E',\phi')$.

The analogous notions are used over a rigid analytic base~$S$.

\medskip
For the following construction we fix a finitely generated projective $A$-submodule $L\subset F^r$ of rank~$r$. Recall that elements of $F^r$ are viewed as row vectors and points in $\Omega^r$ as column vectors. Any $\omega\in\Omega^r$ thus determines an $A$-lattice $L\omega\subset\Cinf$ of rank~$r$. Let $e_{L\omega}$ be the associated exponential function from (\ref{ExpDef}). For any $a\in A\setminus\{0\}$ we have an inclusion of $A$-lattices $L\omega \subset a^{-1}L\omega$ of finite index, so $e_{L\omega}(a^{-1}L\omega)$ is a finite $\BF_q$-subspace of~$\Cinf$. Thus 
\UseTheoremCounterForNextEquation
\begin{equation}\label{PsiLaDef}
\psi^{L\omega}_a\ :=\ a\cdot e_{e_{L\omega}(a^{-1}L\omega)}
\end{equation}
is a polynomial in $\End_{\BF_q}(\BGacomma{\Cinf})$ which by Proposition \ref{exp2} (a) and (b) satisfies the functional equation $\psi^{L\omega}_a (e_{L\omega}(z)) = e_{L\omega}(az)$. Setting also $\psi^{L\omega}_0:=0$, we obtain the Drinfeld $A$-module $(\BGacomma{\Cinf},\psi^{L\omega})$ over $\Cinf$ that is uniformised by the lattice $L\omega$.
As $\omega$ varies over~$\Omega^r$, the exponential function $e_{L\omega}(z)$ is holomorphic in $(z,\omega)\in\Cinf\times\Omega^r$; hence \smash{$\psi^{L\omega}_a$} is holomorphic in $\omega\in\Omega^r$ for each $a\in A$. Together this therefore defines a Drinfeld $A$-module 
\UseTheoremCounterForNextEquation
\begin{equation}\label{PsiLDef}
(\BGacomma{\Omega^r},\psi^{L})
\end{equation}
of rank $r$ over~$\Omega^r$. 

Also, any element $\ell\in F^r$ determines a holomorphic section
\UseTheoremCounterForNextEquation
\begin{equation}\label{PsiLTorsDef}
\mu^L_\ell:\ \omega\mapsto e_{L\omega}(\ell\omega)
\end{equation}
of $\BGacomma{\Omega^r}$ which depends only on the residue class $\ell+L$. 
For any non-zero ideal $N\subset A$ with $N\ell\subset L$ this section lies in the $N$-torsion subgroup $\psi^{L}[N]$ of $\psi^{L}$. Varying $\ell+L$ over $N^{-1}L/L$ this endows the Drinfeld $A$-module $(\BGacomma{\Omega^r},\psi^L)$ with a full level structure of level~$N$ by mapping
\UseTheoremCounterForNextEquation
\begin{equation}\label{PsiLLevelDef}
N^{-1}L/L \longto \psi^L[N],\ \ \ell+L\mapsto \mu^L_\ell.
\end{equation}

Next consider an arbitrary element $\gamma\in\GL_r(F)$. Then for any $\omega\in\Omega^r$ we have $L\omega = L\gamma^{-1}\gamma\omega = j(\gamma,\omega)\cdot L\gamma^{-1}\cdot\gamma(\omega)$ by (\ref{DefOfActionOnOmega}). Multiplication by $j(\gamma,\omega)^{-1}$ thus induces an isomorphism of Drinfeld $A$-modules
\UseTheoremCounterForNextEquation
\begin{equation}\label{PsiLGamma1}
(\BGacomma{\Cinf},\psi^{L\omega})
\ \stackrel{\sim}{\longto}\ 
(\BGacomma{\Cinf},\psi^{L\gamma^{-1}\cdot\gamma(\omega)}).
\end{equation}
Here the target is the pullback of the Drinfeld $A$-module $(\BGacomma{\Omega^r},\psi^{L\gamma^{-1}})$ via the isomorphism $\gamma: \Omega^r\to\Omega^r$, $\omega\mapsto\gamma(\omega)$, evaluated at~$\omega$. Multiplication by the holomorphic function $j(\gamma,{\underline{\ \ }})^{-1}$ thus induces an isomorphism of Drinfeld $A$-modules
\UseTheoremCounterForNextEquation
\begin{equation}\label{PsiLGamma}
(\BGacomma{\Omega^r},\psi^{L})
\ \stackrel{\sim}{\longto}\ 
\gamma^*(\BGacomma{\Omega^r},\psi^{L\gamma^{-1}})
\end{equation}
over~$\Omega^r$. Also, for any $\ell\in F^r$, using Proposition \ref{exp2} (b) we can calculate
\UseTheoremCounterForNextEquation
\begin{equation}\label{PsiLMu}
\begin{array}{rl}
\mu^L_\ell(\omega)
\!\!&=\ e_{L\omega}(\ell\omega) \\
\!\!&=\ e_{j(\gamma,\omega)\cdot L\gamma^{-1}\cdot\gamma(\omega)} \bigl(j(\gamma,\omega)\cdot\ell\gamma^{-1}\cdot\gamma(\omega)\bigr) \\
\!\!&=\ j(\gamma,\omega)\cdot e_{L\gamma^{-1}\cdot\gamma(\omega)} \bigl(\ell\gamma^{-1}\cdot\gamma(\omega)\bigr) \\
\!\!&=\ j(\gamma,\omega)\cdot 
\mu^{L\gamma^{-1}}_{\ell\gamma^{-1}}(\gamma(\omega)).
\end{array}
\end{equation}
Multiplication by $j(\gamma,{\underline{\ \ }})^{-1}$ thus also sends 
the level $N$ structure $\ell+L\mapsto\mu^L_\ell$ of $(\BGacomma{\Omega^r},\psi^{L})$ to the level $N$ structure $\ell\gamma^{-1}+L\gamma^{-1}\mapsto \gamma^*\mu^{L\gamma^{-1}}_{\ell\gamma^{-1}}$ of $\gamma^*(\BGacomma{\Omega^r},\psi^{L\gamma^{-1}})$.

\medskip
Now let $\Gamma<\GL_r(F)$ be an arithmetic subgroup whose right action on $F^r$ normalises the lattice~$L$. Recall from \cite[\PropositionCite 6.2]{Drinfeld1} that $\Gamma<\GL_r(F)$ acts discontinuously on~$\Omega^r$; hence the quotient $\Gamma\backslash \Omega^r$ exists as a rigid analytic space by \cite[\S6.4]{FresnelvdPut}. 
Let $\pi_\Gamma: \Omega^r\onto\Gamma\backslash\Omega^r$ denote the projection morphism.

Assume that $\Gamma$ acts freely on~$\Omega^r$. Then $\Gamma$ also acts freely on $\BGacomma{\Omega^r} = \BGa\times\Omega^r$ through $\gamma(z,\omega) := (j(\gamma,\omega)^{-1}z,\gamma(\omega))$, so the quotient $E_\Gamma := \Gamma\backslash(\BGa\times\Omega^r)$ exists and is a line bundle on $\Gamma\backslash\Omega^r$. By construction
the space of its sections over any open subset $U\subset\Gamma\backslash\Omega^r$ is
\UseTheoremCounterForNextEquation
\begin{equation}\label{ELDef}
E_\Gamma(U)\ :=\ \bigl\{ f:\ \pi_\Gamma^{-1}(U)\to\Cinf\text{\ holomorphic\ }
\bigm| \forall\gamma\in\Gamma:\ f(\gamma(\omega))=j(\gamma,\omega)^{-1} f(\omega) \bigr\}.
\end{equation}
This line bundle comes with a natural isomorphism 
\UseTheoremCounterForNextEquation
\begin{equation}\label{ELPB}
\pi_\Gamma^*E_\Gamma\ \stackrel{\sim}{\longto}\ \BGacomma{\Omega^r}.
\end{equation}
For any $\gamma\in\Gamma$ the equality $\pi_\Gamma=\pi_\Gamma\circ\gamma$ induces a commutative diagram
\UseTheoremCounterForNextEquation
\begin{equation}\label{VoteForMe}
\vcenter{\xymatrix@C-8pt{
\ \pi_\Gamma^*E_\Gamma\ \ar[rrrr]_-\sim^-{(\ref{ELPB})} \ar@{=}[d]
&&&& \ \BGacomma{\Omega^r}\ \ar[d]^\wr \\
\ \gamma^*\pi_\Gamma^*E_\Gamma\ \ar[rrr]_-\sim^-{(\ref{ELPB})}
&&& \ \gamma^*\BGacomma{\Omega^r}\ \ar@{=}[r]
& \ \BGacomma{\Omega^r}\rlap{,}\ \\}}
\end{equation}
where the vertical map on the right is multiplication by $j(\gamma,\underline{\ \ })^{-1}$.
The isomorphism (\ref{PsiLGamma}) for all $\gamma\in\Gamma$ implies that there is a unique Drinfeld $A$-module of the form $(E_\Gamma,\bar\psi^{L})$ over $\Gamma\backslash\Omega^r$ such that (\ref{ELPB}) induces an isomorphism
\UseTheoremCounterForNextEquation
\begin{equation}\label{PiLDrin}
\pi_\Gamma^*(E_\Gamma,\bar\psi^{L}) 
\ \stackrel{\sim}{\longto}\ (\BGacomma{\Omega^r},\psi^{L}).
\end{equation}

Moreover, since $\Gamma$ normalises $L$, it acts on $N^{-1}L/L$ for any non-zero ideal $N\subset A$. For any residue class $\ell+L$ that is fixed by~$\Gamma$, the formula (\ref{PsiLMu}) implies that the associated torsion point $\mu^{L}_\ell$ descends to a torsion point $\bar\mu^L_\ell$ of $(E_\Gamma,\bar\psi^{L})$. In particular, if $\Gamma$ acts trivially on $N^{-1}L/L$, the level $N$ structure (\ref{PsiLLevelDef}) descends to a unique level $N$ structure of $(E_\Gamma,\bar\psi^{L})$
\UseTheoremCounterForNextEquation
\begin{equation}\label{BarPsiLLevelDef}
N^{-1}L/L \longto \bar\psi^L[N],\ \ \ell+L\mapsto\bar\mu^L_\ell.
\end{equation}

\medskip
Now set $\Gamma_U := \Gamma\cap U(F)$ as in (\ref{DefGammaU}) and let $\Lambda' := \iota^{-1}(\Gamma_U) \subset F^{r-1}$ be the corresponding subgroup from (\ref{DefLambdaPrime}), which is commensurable with $A^{r-1}$. Then by Theorem \ref{Thm:Parameter} there exist an admissible open subset $\CU\subset\Cinf\times\Omega^{r-1}$ containing $\{0\}\times\Omega^{r-1}$ and a holomorphic map
\UseTheoremCounterForNextEquation
\begin{equation}\label{Thm:Parameter:Again}
\theta:\ \Gamma_U\backslash\Omega^r \longto \CU,\ \ 
\textstyle \left[\binom{\omega_1}{\omega'}\right] 
\longmapsto \binom{e_{\Lambda'\omega'}(\omega_1)^{-1}}{\omega'}
\end{equation}
which induces an isomorphism of rigid analytic spaces $\Gamma_U\backslash\Omega^r
\stackrel{\sim}{\longto}\CU\cap(\Cinf^\times\times\Omega^{r-1})$.
Also $\pi_\Gamma$ factors through projection morphisms
$$\xymatrix@C+10pt{
\Omega^r \ar[r]^-{\pi_{\Gamma_U}} \ar@/_15pt/[rr]_-{\pi_\Gamma}
& \Gamma_U\backslash\Omega^r  \ar[r]^-{\pi_\Gamma^{\Gamma_U}}
& \Gamma\backslash\Omega^r. \\}$$
For all $\gamma\in\Gamma_U$, the definition (\ref{DefOfj}) implies that $j(\gamma,\omega)=1$ and hence $e_{L\gamma(\omega)} = e_{L\omega}$ and $\psi_a^{L\gamma(\omega)} = \psi_a^{L\omega}$. For ease of notation we denote the function on $\BGa\times\Gamma_U\backslash\Omega^r$ induced by $\psi_a^{L\omega}$ again by $\psi_a^{L\omega}$. Then the Drinfeld $A$-module $(\BGacomma{\Omega^r},\psi^{L})$ is the pullback under $\pi_{\Gamma_U}$ of a unique Drinfeld $A$-module of the form $(\BGacomma{\Gamma_U\backslash\Omega^r},\psi^{L})$ over $\Gamma_U\backslash\Omega^r$. Moreover the isomorphism (\ref{PiLDrin}) descends to a natural isomorphism
\UseTheoremCounterForNextEquation
\begin{equation}\label{PiLDrinU}
(\pi_\Gamma^{\Gamma_U})^*(E_\Gamma,\bar\psi^{L}) 
\ \stackrel{\sim}{\longto}\ (\BGacomma{\Gamma_U\backslash\Omega^r},\psi^{L}).
\end{equation}

\begin{Prop}\label{GenDrinModCodim1}
There exists a unique generalised Drinfeld $A$-module of the form $(\BGacomma{\CU},\tilde\psi^{L})$ over $\CU$ such that 
$$(\BGacomma{\Gamma_U\backslash\Omega^r},\psi^{L})
\ =\ \theta^*(\BGacomma{\CU},\tilde\psi^{L}).$$
Its restriction to $\{0\}\times\Omega^{r-1} \subset \CU$ is a Drinfeld $A$-module of constant rank~$r-1$.
\end{Prop}

\begin{Proof}
We will show that the exponential function ${\Cinf\times(\CU\cap(\Cinf^\times\times\Omega^{r-1})) \longto\Cinf}$, $(z,\theta([\omega]))\mapsto e_{L\omega}(z)$ associated to the Drinfeld $A$-module extends to a holomorphic function on $\Cinf\times\CU$. Writing $\omega=\binom{\omega_1}{\omega'}$ as before, we will express $e_{L\omega}(z)$ as an infinite product in the variables $(z,u,\omega')$ for $u = u_{\omega'}(\omega_1) := e_{\Lambda'\omega'}(\omega_1)^{-1}$ and show that this product also converges near $u=0$. 

For this we define subgroups $L'$ and $L_1$ by the commutative diagram with exact rows
$$\xymatrix@R-24pt{
0 \ar[r] & F^{r-1} \ar[rr] && F^r \ar[rr] && F \ar[r] & 0\\
& \cup \ar@{}[rr]|{v'\mapsto(0,v')} 
&& \cup \ar@{}[rr]|{(v_1,v')\mapsto v_1} && \cup & \\
0 \ar[r] & L' \ar[rr] && L \ar[rr] && L_1 \ar[r] & 0\rlap{.}\\}$$
Since $L$ is commensurable with $A^r$, the subgroups $L'$ and $L_1$ are commensurable with $A^{r-1}$ and $A$, respectively. 
Next, for any $(\ell_1,v')\in L$ and any $\lambda'\in\Lambda'$ we have $\binom{1\ \lambda'}{0\kern5pt 1\kern2pt}\in\Gamma_U$ and hence $(\ell_1,v') \binom{1\ \lambda'}{0\kern5pt 1\kern2pt} = (\ell_1,\ell_1\lambda'+v')\in L$. In particular this implies that $\ell_1\Lambda'\subset L'$. As both $\Lambda'$ and $L'$ are commensurable with $A^{r-1}$, this is an inclusion of finite index if $\ell_1\not=0$.

%
%

Next we fix a subgroup $\tilde L_1\subset L$ which maps isomorphically to $L_1$ under the projection $F^r\onto F$. Then for any $\omega=\binom{\omega_1}{\omega'}\in\Omega^r$ we have $L\omega = \tilde L_1\omega\oplus L'\omega'$. 
Using Proposition \ref{exp2} (a) and the definition (\ref{ExpDef}) of the exponential function, for any $z\in\Cinf$ we thus have 
\UseTheoremCounterForNextEquation
\begin{equation}\label{GenDrinModCodim1Gna1}
e_{L\omega}(z) \ =\ e_{e_{L'\omega'}(L\omega)}(e_{L'\omega'}(z))
\ =\ \tilde z \cdot \!\!\!\!\!\prod_{\ell\in\tilde L_1\setminus\{0\}}
\!\Bigl(1-\frac{\tilde z}{e_{L'\omega'}(\ell\omega)}\Bigr)
\end{equation}
with $\tilde z = e_{L'\omega'}(z)$.
To transform the denominator write $\ell \in \tilde L_1\setminus\{0\}$ in the form $\ell=(\ell_1,v')$ with $\ell_1\in L_1\setminus\{0\}$ and $v'\in F^{r-1}$.
Then we have an inclusion of lattices $\Lambda'\omega' \subset \ell_1^{-1}L'\omega'$, and by the $\Finf$-linear independence of the coefficients of $\omega'$ the index is precisely ${[L':\ell_1\Lambda']}<\infty$.
By the additivity of the exponential function we have
$$e_{\Lambda'\omega'}(\ell_1^{-1}\ell\omega)
\ =\ e_{\Lambda'\omega'}(\omega_1+\ell_1^{-1}v'\omega')
\ =\ u^{-1}+e_{\Lambda'\omega'}(\ell_1^{-1}v'\omega')$$
with $u=e_{\Lambda'\omega'}(\omega_1)^{-1}$. 
Using Proposition \ref{exp2} again we deduce that
\begin{eqnarray*}
e_{L'\omega'}(\ell\omega)
&\!\!=\!\!& \ell_1\cdot e_{\ell_1^{-1}L'\omega'}(\ell_1^{-1}\ell\omega) \\
&\!\!=\!\!& \ell_1\cdot
e_{e_{\Lambda'\omega'}(\ell_1^{-1}L'\omega')}\bigl(
e_{\Lambda'\omega'}(\ell_1^{-1}\ell\omega)\bigr) \\
&\!\!=\!\!& \ell_1\cdot
e_{e_{\Lambda'\omega'}(\ell_1^{-1}L'\omega')}\bigl(
u^{-1}+e_{\Lambda'\omega'}(\ell_1^{-1}v'\omega')\bigr).
\end{eqnarray*}
By the definition and the additivity of the exponential function this in turn yields
\begin{eqnarray*}
e_{L'\omega'}(\ell\omega)
&\!\!=\!\!& \ell_1\cdot
\bigl(u^{-1}+e_{\Lambda'\omega'}(\ell_1^{-1}v'\omega')\bigr) \cdot\kern-10pt
\prod_{{\ell'\in L'\setminus\ell_1\Lambda'}\atop{{\rm modulo}\ \ell_1\Lambda'}}\!\!
\left(1 - \frac{u^{-1}+e_{\Lambda'\omega'}(\ell_1^{-1}v'\omega')}
{e_{\Lambda'\omega'}(\ell_1^{-1}\ell'\omega')}\right) \\
&\!\!=\!\!& \ell_1\cdot
\bigl(u^{-1}+e_{\Lambda'\omega'}(\ell_1^{-1}v'\omega')\bigr) \cdot\kern-10pt
\prod_{{\ell'\in L'\setminus\ell_1\Lambda'}\atop{{\rm modulo}\ \ell_1\Lambda'}}\kern-10pt
\frac{e_{\Lambda'\omega'}(\ell_1^{-1}(\ell'-v')\omega')-u^{-1}}
{e_{\Lambda'\omega'}(\ell_1^{-1}\ell'\omega')} \\
&\!\!=\!\!& \ell_1\cdot
\frac{1+e_{\Lambda'\omega'}(\ell_1^{-1}v'\omega')\cdot u}
{u^{[L':\ell_1\Lambda']}}\cdot\kern-10pt
\prod_{{\ell'\in L'\setminus\ell_1\Lambda'}\atop{{\rm modulo}\ \ell_1\Lambda'}}\kern-10pt
\frac{e_{\Lambda'\omega'}(\ell_1^{-1}(\ell'-v')\omega')\cdot u-1}
{e_{\Lambda'\omega'}(\ell_1^{-1}\ell'\omega')} \\
&\!\!=\!\!& \frac{\ell_1}{u^{[L':\ell_1\Lambda']}}\cdot
\frac
{\displaystyle\prod_{{\ell'\in L'}\ {{\rm mod}\ \ell_1\Lambda'}}
\kern-10pt \bigl(1-e_{\Lambda'\omega'}(\ell_1^{-1}(\ell'-v')\omega')\cdot u\bigr)}
{\displaystyle\prod_{{\ell'\in L'\setminus \ell_1\Lambda'}\ {{\rm mod}\ \ell_1\Lambda'}}
\kern-20pt e_{\Lambda'\omega'}(\ell_1^{-1}\ell'\omega')\kern30pt},
\end{eqnarray*}
where the last transformation uses the fact that $(-1)^{[L':\ell_1\Lambda']-1}=1$ because $[L':\ell_1\Lambda']$ is a power of~$q$. 
Plugging this into the formula (\ref{GenDrinModCodim1Gna1}) we conclude that
\UseTheoremCounterForNextEquation
\begin{equation}\label{GenDrinModCodim1Gna2}
e_{L\omega}(z) \ =\ \tilde z \cdot \kern-15pt
\prod_{(\ell_1,v')\in\tilde L_1\setminus\{0\}} \kern-2pt
\left(1-\tilde z\cdot 
\frac{u^{[L':\ell_1\Lambda']}}{\ell_1}\cdot
\frac
{\displaystyle\prod_{{\ell'\in L'\setminus \ell_1\Lambda'}\ {{\rm mod}\ \ell_1\Lambda'}}
\kern-20pt e_{\Lambda'\omega'}(\ell_1^{-1}\ell'\omega')\kern30pt}
{\displaystyle\prod_{{\ell'\in L'}\ {{\rm mod}\ \ell_1\Lambda'}}
\kern-10pt \bigl(1-e_{\Lambda'\omega'}(\ell_1^{-1}(\ell'-v')\omega')\cdot u\bigr)}
\right)_.
\end{equation}
As $(\ell_1,\ell')$ runs through $\tilde L_1\setminus\{0\}$, the index $[L':\ell_1\Lambda']$ goes to infinity. Using the geometric series we can therefore expand the right hand side of (\ref{GenDrinModCodim1Gna2}) as a power series in $u$ whose coefficients are functions of $(\tilde z,\omega_1)$. 
We will show that this expression converges locally uniformly for all $\tilde z\in\Cinf$ and all $(u,\omega_1)$ in a suitable tubular neighbourhood of $\{0\}\times\Omega^{r-1}$.

For this take any $n>0$. By Proposition \ref{NewEarth} (c) there exists a constant $c_n>0$, such that for any $\omega'\in\Omega^{r-1}_n$ and any $v'\in F_\infty^{r-1}$ we have $|e_{\Lambda'\omega'}(v'\omega')|<c_n$. In particular this inequality holds for $\ell_1^{-1}\ell'$ and $\ell_1^{-1}(\ell'-v')$ in place of $v'$. Thus if $|u|\le r_n := (2c_n)^{-1}$, we have $|e_{\Lambda'\omega'}(\ell_1^{-1}(\ell'-v')\omega')\cdot u|<2^{-1}$, so the geometric series for 
$$\frac{1}{1-e_{\Lambda'\omega'}(\ell_1^{-1}(\ell'-v')\omega')\cdot u}$$
converges uniformly to a value of norm~$1$. 
Combining the inequalities yields the bound
$$\left|\frac{u^{[L':\ell_1\Lambda']}}{\ell_1}\cdot
\frac
{\displaystyle\prod_{{\ell'\in L'\setminus \ell_1\Lambda'}\ {{\rm mod}\ \ell_1\Lambda'}}
\kern-20pt e_{\Lambda'\omega'}(\ell_1^{-1}\ell'\omega')\kern30pt}
{\displaystyle\prod_{{\ell'\in L'}\ {{\rm mod}\ \ell_1\Lambda'}}
\kern-10pt \bigl(1-e_{\Lambda'\omega'}(\ell_1^{-1}(\ell'-v')\omega')\cdot u\bigr)}
\right|
\ \le\ \frac{r_n^{[L':\ell_1\Lambda']}c_n^{[L':\ell_1\Lambda']-1}}{|\ell_1|}
\ =\ \frac{2^{-[L':\ell_1\Lambda']}}{|\ell_1|c_n}.$$
As both $|\ell_1|$ and $[L':\ell_1\Lambda']$ go to infinity with~$\ell_1$, for any $R>0$ this proves that the right hand side of (\ref{GenDrinModCodim1Gna2}) converges uniformly for all $(\tilde z,u,\omega') \in B(0,R)\times B(0,r_n)\times\Omega^{r-1}_n$. Varying $n$ and $R$ it therefore converges locally uniformly on $\Cinf\times\CT$ for the tubular neighbourhood $\CT := \bigcup_{n\ge1}B(0,r_n)\times\Omega^{r-1}_n$ and the limit is a holomorphic function of $(\tilde z,u,\omega')$. 
Substituting $\tilde z = e_{L'\omega'}(z)$, which is already a holomorphic function of $(z,\omega')\in\Cinf\times\Omega^{r-1}$, thus yields a holomorphic function $E(z,u,\omega')$ on $\Cinf\times\CT$ such that 
\UseTheoremCounterForNextEquation
\begin{equation}\label{GenDrinModCodim1Gna3}
e_{L\omega}(z)\ =\ E(z,e_{\Lambda'\omega'}(\omega_1)^{-1},\omega')
\end{equation}
for all $z\in\Cinf$ and $\omega=\binom{\omega_1}{\omega'}\in\Omega^r$ with $\theta([\omega])\in\CT$.
Now recall that for any $\omega\in\Omega^r$, the Drinfeld $A$-module $\psi^{L\omega}$ is characterised by the fact that for each $a\in A\setminus\{0\}$ the function $\psi_a^{L\omega}$ is an $\BF_q$-linear polynomial in $\Cinf[z]$ satisfying the functional equation $\psi_a^{L\omega}(e_{L\omega}(z)) = e_{L\omega}(az)$. 
Writing this as an identity of power series in~$z$ and observing that $e_{L\omega}(z)=z+($higher terms), it follows that each coefficient of $\psi_a^{L\omega}$ is a certain polynomial with coefficients in $A$ in finitely many coefficients of $e_{L\omega}(z)$. By what we have just proved, these coefficients, as functions of $(e_{\Lambda'\omega'}(\omega_1)^{-1},\omega')$, extend to holomorphic functions of $(u,\omega')\in\CT$. Thus the same is true for the coefficients of $\psi_a^{L\omega}$. In other words, there is a unique holomorphic function $\tilde\psi_a^{L}$ on $\Cinf\times\CT$, which is an $\BF_q$-linear polynomial of degree $\le r\deg(a)$ in~$z$, such that 
\UseTheoremCounterForNextEquation
\begin{equation}\label{GenDrinModCodim1Gna4}
\psi_a^{L\omega}(z)\ =\ \tilde\psi_a^{L}(z,e_{\Lambda'\omega'}(\omega_1)^{-1},\omega')
\end{equation}
for all $z\in\Cinf$ and $\omega=\binom{\omega_1}{\omega'}\in\Omega^r$ with $\theta([\omega])\in\CT$. Setting $\tilde\psi_0^{L}:=0$, the fact that $a\mapsto\psi_a^{L}$ is an $\BF_q$-algebra homomorphism by continuity implies that $a\mapsto\tilde\psi_a^{L}$ is also $\BF_q$-algebra homomorphism. Moreover, the fact that $\tfrac{\partial}{\partial z}\psi^L_a=a$ identically implies that $\tfrac{\partial}{\partial z}\tilde\psi^L_a=a$ identically as well. Furthermore, by continuity the functional equation $\psi_a^{L\omega}(e_{L\omega}(z)) = e_{L\omega}(az)$ extends to a functional equation 
\UseTheoremCounterForNextEquation
\begin{equation}\label{GenDrinModCodim1Gna5}
\tilde\psi^{L}_a\bigl(E(z,u,\omega'),u,\omega'\bigr)\ =\ E(az,u,\omega')
\end{equation}
for all $z\in\Cinf$ and $(u,\omega')\in\CT$. 
If we substitute $u:=0$, the right hand side of (\ref{GenDrinModCodim1Gna2}) becomes just $\tilde z = e_{L'\omega'}(z)$; hence $E(z,0,\omega')=e_{L'\omega'}(z)$. Thus (\ref{GenDrinModCodim1Gna5}) reduces to the equation
\UseTheoremCounterForNextEquation
\begin{equation}\label{GenDrinModCodim1Gna6}
\tilde\psi^{L}_a\bigl(e_{L'\omega'}(z),0,\omega'\bigr)\ =\ e_{L'\omega'}(az).
\end{equation}
For any $\omega'\in\Omega^{r-1}$ the map $a\mapsto\tilde\psi^{L}_a(\underline{\ \ },0,\omega')$ is therefore the Drinfeld $A$-module of rank $r-1$ associated to the lattice $L'\omega'\subset\Cinf$. 
All this together proves that $a\mapsto\tilde\psi^{L}_a$ constitutes a generalised Drinfeld $A$-module of rank $\le r$ over $\CT$, whose restriction to the locus $u=0$ is a Drinfeld $A$-module of constant rank~$r-1$. 

We have thus proved the desired statement over~$\CT$. Since $\tilde\psi^{L}$ is already given over $\CU\cap(\Cinf^\times\times\Omega^{r-1})$, the existence and uniqueness also follows over~$\CU$, as desired.
\end{Proof}

\section{Drinfeld moduli spaces}
\label{Sec:DMS}

Let $\hatA\cong\smash{\prod_{\Fp}A_\Fp}$ be the profinite completion of~$A$ and $\BAfF = \hatA\otimes_A F$ the ring of finite ad\`eles of~$F$. 
For any open compact subgroup $K<\GL_r(\hatA)$ 
let $M^r_{A,K}$ be the  \emph{Drinfeld modular variety of level~$K$}, which is a normal integral affine algebraic variety over~$F$. 
The associated rigid analytic space over $\Cinf$ possesses a natural isomorphism
\UseTheoremCounterForNextEquation
\begin{equation}\label{AnalUniv1}
\GL_r(F)\!\bigm\backslash\! \bigl( \Omega^r\times\GL_r(\BAfF)/K \bigr)
\ \stackrel{\sim}{\longto}\ M^r_{A,K}(\Cinf),
\end{equation}
whose precise characterisation we shall describe below.
For any $g\in\GL_r(\BAfF)$ let $\pi_g$ denote the composite morphism
\UseTheoremCounterForNextEquation
\begin{equation}\label{AnalUnivPig1b1}
\xymatrix@R-25pt@C+5pt{
\Omega^r \ar[r]
& \GL_r(F)\!\bigm\backslash\! \bigl( \Omega^r\times\GL_r(\BAfF)/K \bigr) 
\ar[r]_-\sim^-{(\ref{AnalUniv1})} & M^r_{A,K}(\Cinf),\\
[\omega] \ar@{|->}[r]  & [(\omega,g)] . \\}
\end{equation}
Consider the arithmetic subgroup 
\UseTheoremCounterForNextEquation
\begin{equation}\label{AnalUnivLion}
\Gamma_g \ :=\ \GL_r(F)\cap gKg^{-1}.
\end{equation}
Then $\pi_g$ factors through an isomorphism $\Gamma_g\backslash\Omega^r\stackrel{\sim}{\longto} M_g(\Cinf)$ for a unique connected component $M_g$ of $M^r_{A,K}\times_{\Spec F}\Spec\Cinf$. In other words we have a commutative diagram
\UseTheoremCounterForNextEquation
\begin{equation}\label{AnalUnivPig1b2}
\vcenter{\xymatrix@R-5pt@C+20pt{
\Omega^r \ar[r]^-{\pi_g} \ar@{->>}[d]^{\pi_{\Gamma_g}}
& M^r_{A,K}(\Cinf) \\
\Gamma_g\backslash\Omega^r \ar[r]^-\sim \ar@{^{ (}->}[ur]^{i_g}
& M_g(\Cinf). \ar@{}[u]|{\textstyle\cup} \\}}
\end{equation}
For any $\gamma\in\GL_r(F)$ and $k\in K$ we have $[(\omega,g)] = [(\gamma(\omega),\gamma gk)]$ and hence
\UseTheoremCounterForNextEquation
\begin{equation}\label{AnalUnivPig2}
\pi_g = \pi_{\gamma gk}\circ\gamma.
\end{equation}
For any two elements $g$, $g'\in\GL_r(\BAfF)$ we have $M_g=M_{g'}$ if and only if $g$ and $g'$ represent the same double coset in $\GL_r(F)\backslash\GL_r(\BAfF)/K$.
Thus for any choice of representatives $g_1,\ldots,g_n\in\GL_r(\BAfF)$ we have
\UseTheoremCounterForNextEquation
\begin{equation}\label{AnalUnivDec}
M^r_{A,K}\times_{\Spec F}\Spec\Cinf\ =\ \smash{\coprod_{i=1}^n M_{g_i}}.
\end{equation}
Since $M^r_{A,K}$ is integral, these connected components over $\Cinf$ are Galois conjugate over~$F$. Let $F_K$ denote the field of constants of $M^r_{A,K}$ (which is a certain ray class field of~$F$ that can be characterised uniquely by abelian class field theory). Then the different connected components $M_{g_i}$ are just the varieties obtained by base change $M^r_{A,K}\times_{\Spec F_K}\Spec\Cinf$ for all $F$-linear embeddings $F_K\into\Cinf$.

\medskip
For later use we also record:

\begin{Prop}\label{StrongApprox}
Elements $g_1,\ldots,g_n\in\GL_r(\BAfF)$ form representatives of the double quotient $\GL_r(F)\backslash\GL_r(\BAfF)/K$ if and only if their determinants $\det(g_1),\ldots,\det(g_n)$ form representatives of $F^\times\backslash(\BAfF)^\times/\det(K)$.
\end{Prop}

\begin{Proof}
Direct consequence of strong approximation for the simply connected reductive group $\SL_r$ to the effect that the closure of $\SL_r(F)$ in $\GL_r(\BAfF)$ is $\SL_r(\BAfF)$.
\end{Proof}


\medskip
Now assume that $K$ is \emph{fine}, which by \cite[\DefinitionCite1.4]{Pink} means that 
the image of $K$ in $\GL_r(A/\Fp)$ is unipotent for some maximal ideal $\Fp\subset A$. Then by \cite[\PropositionCite1.5]{Pink} there is a natural \emph{universal family of Drinfeld $A$-modules} $(E,\phi)$ over \smash{$M^r_{A,K}$}, using which one can interpret \smash{$M^r_{A,K}$} as a fine moduli space of Drinfeld $A$-modules with some generalised level structure. 
The pullback of $(E,\phi)$ under the morphism (\ref{AnalUniv1}) can be described as follows. Viewing elements of $F^r$ and $\hatA^r$ and $(\BAfF)^r_{\vphantom{t}}$ as row vectors, for any $g\in \GL_r(\BAfF)$ we set 
\UseTheoremCounterForNextEquation
\begin{equation}\label{LatticeDef}
L_g\ :=\ \hatA^rg^{-1}\cap F^r\ \subset\ (\BAfF)^r_{\vphantom{t}},
\end{equation}
which is a finitely generated projective $A$-module of rank~$r$. Since $K<\GL_r(\hatA)$, by construction the right action of $\Gamma_g$ on $F^r$ normalises~$L_g$. Moreover, the assumption that $K$ is fine implies that all torsion elements of $\Gamma_g$ are unipotent; hence $\Gamma_g$ acts freely on $\Omega^r$. There is therefore a natural Drinfeld $A$-module $(E_{\Gamma_g},\bar\psi^{L_g})$ over $\Gamma_g\backslash\Omega^r$ such that $\pi_{\Gamma_g}^*(E_{\Gamma_g},\bar\psi^{L_g})
\cong (\BGacomma{\Omega^r},\psi^{L_g})$ by (\ref{PiLDrin}). For this there is a natural isomorphism
\UseTheoremCounterForNextEquation
\begin{equation}\label{StoneTrek}
i_g^*(E,\phi) 
\ \stackrel{\sim}{\longto}\ (E_{\Gamma_g},\bar\psi^{L_g}).
\end{equation}

Moreover, suppose that $K$ is the principal congruence subgroup of level~$N$
$$K(N)\ :=\ \bigl\{k\in\GL_r(\hatA)\bigm| k\equiv\Id_r\mathrel{\rm mod} N\bigr\}$$
for some non-zero ideal $N\subset A$. Then \smash{$M^r_{A,K(N)}$} represents the functor which to any scheme $S$ over $F$ associates the set of isomorphism classes of tuples $(E,\phi,\mu)$ consisting of a Drinfeld $A$-module $(E,\phi)$ of rank $r$ over $S$ and a full level $N$ structure $\mu: N^{-1}A^r/A^r \to \phi[N]$. For any $g\in\GL_r(\BAfF)$ we then have 
$$\Gamma_g\ =\ \bigl\{\gamma\in\GL_r(F)\bigm| (\ell+L_g)\gamma=\ell+L_g\ \hbox{for all\ }\ell\in N^{-1}L_g\bigr\}.$$
Thus the Drinfeld $A$-module $(E_{\Gamma_g},\bar\psi^{L_g})$ on $\Gamma_g\backslash\Omega^r$ is endowed with a full level $N$ structure $\bar\mu^{L_g}: N^{-1}L_g/L_g \to \bar\psi^{L_g}[N]$ by (\ref{BarPsiLLevelDef}).
To any coset $\ell+A^r \subset N^{-1}A^r$ associate the coset
\UseTheoremCounterForNextEquation
\begin{equation}\label{Aurora1}
\ell_g+L_g\ :=\ (\ell+\hatA^r)g^{-1}\cap F^r
\ \subset\ N^{-1}L_g.
\end{equation}
This induces an isomorphism $N^{-1}A^r/A^r \stackrel{\sim}{\to} N^{-1}L_g/L_g$. 
The isomorphism (\ref{StoneTrek}) sends the level $N$ structure $\ell+A^r\mapsto i_g^*\mu(\ell+A^r)$ to the level $N$ structure $\ell+A^r \mapsto \ell_g+L_g\mapsto \bar\mu^L_\ell$. In fact this characterises the isomorphism (\ref{StoneTrek}) uniquely. Moreover, since \smash{$M^r_{A,K(N)}$} is a fine moduli space for Drinfeld $A$-modules with a full level $N$ structure, this also characterises the isomorphism (\ref{AnalUniv1}) uniquely in this case.

\medskip
For an arbitrary open compact subgroup~$K$, choose any $N$ such that ${K(N)\triangleleft K}$. Then the finite group $K/K(N)$ acts on $M^r_{A,K(N)}$ by transforming the level $N$ structure, and the quotient is naturally isomorphic to $M^r_{A,K}$. The group $K/K(N)$ also acts by right multiplication on $\GL_r(F)\backslash(\Omega^r\times\GL_r(\BAfF)/K(N))$, and the isomorphism (\ref{AnalUniv1}) in the case of $K$ is obtained from that in the case of $K(N)$ by taking quotients. In particular, the two instances of the map $i_g$ from (\ref{AnalUnivPig1b2}) for $K$ and $K(N)$ form a commutative diagram with the projection $M^r_{A,K(N)} \onto M^r_{A,K}$. 

Similarly, if $K$ is fine, in \cite[\PropositionCite 1.5]{Pink} the universal family on $M^r_{A,K}$ was constructed precisely so that its pullback is the given universal family over $M^r_{A,K(N)}$. The isomorphism (\ref{StoneTrek}) in the case of $K$ is the unique one whose pullback yields the isomorphism (\ref{StoneTrek}) in the case of $K(N)$.

It is useful to know that isomorphisms of Drinfeld modules can be characterised uniquely by using just one torsion point. Since $K$ is fine, its image in $\GL_r(A/\Fp)$ is unipotent for some maximal ideal $\Fp\subset A$, and so it fixes some non-zero coset $\ell+\hatA^r \subset \Fp^{-1}\hatA^r$. For each 
$g\in\GL_r(\BAfF)$ the subgroup $\Gamma_g$ then fixes the corresponding coset $\ell_g+L_g \subset \Fp^{-1}L_g$ defined by (\ref{Aurora1}). The associated torsion point $\mu^{L_g}_{\ell_g}$ thus descends to a nowhere zero $\Fp$-torsion point of $(E_{\Gamma_g},\bar\psi^{L_g})$ over $\Gamma_g\backslash\Omega^r$. 
On the other hand, choosing $N\subset\Fp$, the group $K/K(N)$ fixes the coset $\ell+\hatA^r$; hence the associated $\Fp$-torsion point coming from the level $N$ structure descends to a nowhere zero $\Fp$-torsion point of the universal family $(E,\phi)$ over \smash{$M^r_{A,K}$}. By construction the isomorphism (\ref{StoneTrek}) identifies the respective $\Fp$-torsion points. As any isomorphism of Drinfeld modules is scalar and hence determined by the image of any non-zero point, it follows that the isomorphism is uniquely characterised by this.

\medskip
In the following we care mostly about the composite isomorphism 
\UseTheoremCounterForNextEquation
\begin{equation}\label{Tuvok}
\qquad\qquad\xymatrix@C+20pt{
\llap{$\pi_g^*(E,\phi)\ =\ $}
\pi_{\Gamma_g}^* i_g^* (E,\phi) \ar[r]_-\sim^-{(\ref{StoneTrek})} 
& \pi_{\Gamma_g}^* (E_{\Gamma_g},\bar\psi^{L_g}) \ar[r]_-\sim^-{(\ref{PiLDrin})} 
& (\BGacomma{\Omega^r},\psi^{L_g}). \\}
\end{equation}
This changes with $g$ as follows. Consider any $g\in\GL_r(\BAfF)$ and $\gamma\in\GL_r(F)$ and $k\in K$. Since $K<\GL_r(\hatA)$, from (\ref{LatticeDef}) we deduce that
$$L_{\gamma gk}
\ =\ \hatA^rk^{-1}g^{-1}\gamma^{-1}\cap F^r
\ =\ (\hatA^rg^{-1}\cap F^r)\gamma^{-1}
\ =\ L_g\gamma^{-1}.$$
The isomorphisms from (\ref{Tuvok}) for $g$ and for $\gamma gk$ thus fit into a diagram
\UseTheoremCounterForNextEquation
\begin{equation}\label{AnalUnivPullback2a}
\vcenter{\xymatrix@C+50pt{
\pi_g^* (E,\phi) \ar[r]_-\sim^-{(\ref{Tuvok})\;{\rm for}\;g} 
\ar@{=}[d]^{(\ref{AnalUnivPig2})}
& (\BGacomma{\Omega^r},\psi^{L_g}) \ar[d]_\wr^{(\ref{PsiLGamma})} \\
\gamma^* \pi_{\gamma gk}^* (E,\phi) \ar[r]_-\sim^-{(\ref{Tuvok})\;{\rm for}\;\gamma gk} 
& \gamma^* (\BGacomma{\Omega^r},\psi^{L_{\gamma gk}}), \\}}
\end{equation}
where the vertical map on the right is multiplication by $j(\gamma,\underline{\ \ })^{-1}$. Using (\ref{PsiLMu}) one verifies that the isomorphisms preserve some nowhere vanishing torsion point. Thus the two composites must coincide; in other words the diagram (\ref{AnalUnivPullback2a}) commutes.

\medskip
We end this section by looking at functoriality. Consider a second open compact subgroup $K'<\GL_r(\hatA)$ and an element $h\in\GL_r(\BAfF)$ such that $hK'h^{-1}<K$. Then there is a well-defined map
\UseTheoremCounterForNextEquation
\begin{equation}\label{Jh1}
\qquad\xymatrix@R-25pt{
\llap{$J_h:\ $}
\GL_r(F)\!\bigm\backslash\! \bigl( \Omega^r\times\GL_r(\BAfF)/K'\bigr)
\ar[r] & 
\GL_r(F)\!\bigm\backslash\! \bigl( \Omega^r\times\GL_r(\BAfF)/K\bigr), \\
[(\omega,gh)] \ar@{|->}[r] & [(\omega,g)]. \\}
\end{equation}
If $h$ has coefficients in~$\hatA$, we have $\hatA^r\subset\hatA^rh^{-1}$ and hence 
$$L_g\ =\ \hatA^rg^{-1}\cap F^r\ \subset\ \hatA^rh^{-1}g^{-1}\cap F^r\ =\ L_{gh}$$
for any $g\in\GL_r(\BAfF)$. Thus for any $\omega\in\Omega^r$ we have $L_g\cdot\omega\subset L_{gh}\cdot\omega$, and using Proposition \ref{exp2} (a) we obtain an isogeny of Drinfeld modules 
\UseTheoremCounterForNextEquation
\begin{equation}\label{xih1}
\tilde\eta_h\ :=\ e_{e_{L_g\cdot\omega}(L_{gh}\cdot\omega)}:\ 
(\BGacomma{\Omega^r},\psi^{L_g}) \longto (\BGacomma{\Omega^r},\psi^{L_{gh}}).
\end{equation}
By contrast, if $h^{-1}$ has coefficients in~$\hatA$, we have $\hatA^rh^{-1}\subset\hatA^r$ and hence $L_{gh}\subset L_g$, which yields an isogeny of Drinfeld modules 
\UseTheoremCounterForNextEquation
\begin{equation}\label{xih2}
\tilde\xi_h\ :=\ e_{e_{L_{gh}\cdot\omega}(L_g\cdot\omega)}:\ 
(\BGacomma{\Omega^r},\psi^{L_{gh}}) \longto (\BGacomma{\Omega^r},\psi^{L_g}).
\end{equation}
By construction the isogenies $\tilde\eta_h$ and $\tilde\xi_h$ are mutually inverse isomorphisms if $h\in\GL_r(\hatA)$. In analogy with (\ref{AnalUnivPig1b1}) write 
$$\qquad\xymatrix@R-25pt@C+5pt{
\llap{$\pi'_{gh}:\ \ $}\Omega^r \ar[r]
& \GL_r(F)\!\bigm\backslash\! \bigl( \Omega^r\times\GL_r(\BAfF)/K' \bigr) 
\ar[r]_-\sim^-{\smash{(\ref{AnalUniv1})}} & M^r_{A,K'}(\Cinf),\\
[\omega] \ar@{|->}[r]  & [(\omega,gh)] . \\}$$

\begin{Prop}\label{Func1}
\begin{itemize}
\item[(a)] Via (\ref{AnalUniv1}) the map $J_h$ corresponds to a morphism of varieties
$$\smash{J_h:\ M^r_{A,K'}\longto M^r_{A,K}.}$$
\item[(b)] For every $g\in\GL_r(\BAfF)$ we have $\pi_g = J_h\circ \pi'_{gh}$.
\end{itemize}
\medskip\noindent Now assume that $K$ and $K'$ are fine, and let $(E,\phi)$ and $(E',\phi')$ denote the respective universal families on $M^r_{A,K}$ and $M^r_{A,K'}$. Then:
\begin{itemize}
\item[(c)] If $h$ has coefficients in~$\hatA$, there is a natural isogeny $\eta_h: J_h^*(E,\phi) \to (E',\phi')$ which for every $g\in\GL_r(\BAfF)$ makes the following diagram commute:
$$\xymatrix{
\pi_g^*(E,\phi) \ar@{=}[r]^-{\rm(b)} \ar[d]_\wr^{(\ref{Tuvok})\;{\rm for}\;g} &
\pi_{gh}^{\prime*}J_h^*(E,\phi) 
\ar[rr]^{\pi_{gh}^{\prime*}\eta_h} &&
\pi_{gh}^{\prime*}(E',\phi') \ar[d]_\wr^{(\ref{Tuvok})\;{\rm for}\;gh} \\
(\BGacomma{\Omega^r},\psi^{L_g}) \ar[rrr]^{\tilde\eta_h} &&&
(\BGacomma{\Omega^r},\psi^{L_{gh}}).\\}$$
\item[(d)] If $h^{-1}$ has coefficients in~$\hatA$, there is a natural isogeny $\xi_h: (E',\phi')\to J_h^*(E,\phi)$ which for every $g\in\GL_r(\BAfF)$ makes the following diagram commute:
$$\xymatrix{
\pi_{gh}^{\prime*}(E',\phi') \ar[d]_\wr^{(\ref{Tuvok})\;{\rm for}\;gh} 
\ar[rr]^{\pi_{gh}^{\prime*}\xi_h} &&
\pi_{gh}^{\prime*}J_h^*(E,\phi) \ar@{=}[r]^-{\rm(b)} &
\pi_g^*(E,\phi) \ar[d]_\wr^{(\ref{Tuvok})\;{\rm for}\;g} \\
(\BGacomma{\Omega^r},\psi^{L_{gh}}) \ar[rrr]^{\tilde\xi_h} &&& 
(\BGacomma{\Omega^r},\psi^{L_g}) .\\}$$
\item[(e)] For any $a\in A\setminus\{0\}$ such that both $h$ and $ah^{-1}$ have coefficients in~$\hatA$, we have $\eta_h\circ\xi_{a^{-1}h}=\phi'_a$ and $\xi_{a^{-1}h}\circ\eta_h = J_h^*\phi_a$.
\item[(f)] If $h\in\GL_r(\BAfF)$ is a scalar matrix and $K=K'$, then $J_h$ is the identity morphism. If in addition $h=a\cdot\Id_r$ for $a\in A\setminus\{0\}$, then $\eta_h=\phi_a$. If instead $h=a^{-1}\cdot\Id_r$ for $a\in A\setminus\{0\}$, then $\xi_h=\phi_a$.
\end{itemize}
\end{Prop}

\begin{Proof}
(Sketch) The formulas in (b), (e), and (f) follow by direct calculation from the constructions in (\ref{Jh1}) and (\ref{xih1}) and (\ref{xih2}), once the remaining assertions are proved.

The constructions of $J_h$ and $\xi_h$ in (a) and (d) are those of \cite[\PropositionsCite2.6--7]{Pink}. 
(Except that due to the change of convention explained in Remark \ref{ConventionsRemA} the present morphism $J_h$ corresponds to the morphism $J_{(h^T)^{-1}}$ from \cite[\PropositionCite2.6]{Pink}, and the present isogeny $\xi_h$ to the isogeny $\xi_{(h^T)^{-1}}$ from \cite[\PropositionCite2.7]{Pink}.) Roughly speaking, by taking invariants everything reduces to the case that $K=K(N)$ and $K'=K(N')$, where $J_h$ and $\xi_h$ can be described explicitly using the modular interpretation. 

The construction of $\eta_h$ in (c) is dual to that of $\xi_h$ and follows the same principles. For an alternative construction observe that the formulas in (e) characterise $\eta_h$ uniquely in terms of $\xi_{a^{-1}h}$. Noting that the endomorphism $\phi'_a$ of $(E',\phi')$ also factors through the isogeny $\xi_{a^{-1}h}: (E',\phi')\to J_h^*(E,\phi)$ constructed via the modular interpretation, one can construct $\eta_h$ by the formula $\eta_h\circ\xi_{a^{-1}h}=\phi'_a$ and deduce its properties from that.
\end{Proof}

\begin{Prop}\label{Func2}
Consider open compact subgroups $K,K',K''<\GL_r(\hatA)$ and elements $h,h'\in\GL_r(\BAfF)$ such that $hK'h^{-1}<K$ and $h'K''h^{\prime-1}<K'$. Then we have:
\begin{itemize}
\item[(a)] $J_{hh'}=J_h\circ J_{h'}$.
\item[(b)] $\eta_{hh'} = \eta_{h'}\circ J_{h'}^*\eta_h$ if $K,K',K''$ are fine and $h,h'$ have coefficients in~$\hatA$.
\item[(c)] $\xi_{hh'} = J_{h'}^*\xi_h\circ\xi_{h'}$ if $K,K',K''$ are fine and $h^{-1},h^{\prime-1}$ have coefficients in~$\hatA$.
\end{itemize}
\end{Prop}

\begin{Proof}
Direct calculation for the maps in (\ref{Jh1}) and (\ref{xih1}) and (\ref{xih2}).
\end{Proof}

\section{Satake compactification}
\label{Sec:Satake}

According to \cite[\DefinitionCite 4.1]{Pink}, any normal integral proper algebraic variety $\OM^r_{A,K}$ over $F$ which contains \smash{$M^r_{A,K}$} as an open dense subvariety, such that the universal family $(E,\phi)$ extends to a weakly separating generalised Drinfeld $A$-module $(\bar{E},\bar{\phi})$ over \smash{$\OM^r_{A,K}$}, is called a \emph{Satake compactification of $M^r_{A,K}$.} By \cite[\TheoremCite 4.2]{Pink}, such a Satake compactification exists and is projective over~$F$, and together with its ``universal family'' $(\bar{E},\bar{\phi})$ it is uniquely determined up to unique isomorphism. The proof, however, tells us very little about what the boundary of this compactification looks like. 

A rigid analytic construction of the same Satake compactification was given by Kapranov \cite{Kapranov} in the special case $A=\BF_q[t]$ and by H\"aberli \cite{Haeberli} in general. They explicitly construct a rigid analytic space that is projective over~$\Cinf$ and has a natural stratification by finitely many rigid analytic spaces of the form $\Gamma'\backslash\Omega^{r'}$ for integers $1\le r'\le r$ and arithmetic subgroups $\Gamma'<\GL_{r'}(F)$. H\"aberli also proves that the result is naturally isomorphic to $\OM^r_{A,K}(\Cinf)$. What we need from this is an analytic description of $\OM^r_{A,K}$ along all boundary strata of codimension~$1$, where 
the fibers of the universal family $(\bar{E},\bar{\phi})$ are Drinfeld modules of rank $r-1$. 

\medskip
Since $\OM^r_{A,K}$ is integral and contains $M^r_{A,K}$ as an open dense subvariety, each connected component $M_g$ of $M^r_{A,K}\times_{\Spec F}\Spec\Cinf$ is open and dense in a connected component \smash{$\OM_g$} of \smash{$\OM^r_{A,K}\times_{\Spec F}\Spec\Cinf$}, and the decomposition (\ref{AnalUnivDec}) extends to a decomposition
\UseTheoremCounterForNextEquation
\begin{equation}\label{AnalUniv4a}
\OM^r_{A,K}\times_{\Spec F}\Spec\Cinf\ =\ \smash{\coprod_{i=1}^n \OM_{g_i}}.
\end{equation}
Also, the field of constants of $\OM^r_{A,K}$ is again~$F_K$, and the connected components $\OM_{g_i}$ are just the varieties obtained by base change $\OM^r_{A,K}\times_{\Spec F_K}\Spec\Cinf$ for all $F$-linear embeddings $F_K\into\Cinf$.

Assume that $K$ is fine. Consider any $g\in\GL_r(\BAfF)$, and set $\Gamma_{g,U} := \Gamma_g\cap U(F)$ and $\Lambda'_g := \iota^{-1}(\Gamma_{g,U}) \subset F^{r-1}$ as in (\ref{DefGammaU}) and (\ref{DefLambdaPrime}). By Theorem \ref{Thm:Parameter} there exist an admissible open subset $\CU_g\subset\Cinf\times\Omega^{r-1}$ containing $\{0\}\times\Omega^{r-1}$ and a holomorphic map
\UseTheoremCounterForNextEquation
\begin{equation}\label{Thm:ParameterA:YetAgain}
\textstyle\theta_g:\ 
\Gamma_{g,U}\backslash\Omega^r \longto \CU_g,\ \ 
\left[\binom{\omega_1}{\omega'}\right] 
\longmapsto \binom{e_{\Lambda_g'\omega'}(\omega_1)^{-1}}{\omega'},
\end{equation}
which induces an isomorphism of rigid analytic spaces $\Gamma_{g,U}\backslash\Omega^r \stackrel{\sim}{\longto} \CU_g\cap(\Cinf^\times\times\Omega^{r-1})$.

\begin{Prop}\label{Codim1Prop1}
\begin{itemize}
\item[(a)] There exists a unique morphism of rigid analytic spaces $\bar\pi_g: \CU_g\to\OM_g(\Cinf)$ making the following diagram commute:
$$\xymatrix@C+5pt{
\Omega^r \ar@{->>}[r]^-{\pi_{\Gamma_{g,U}}} \ar[d]^{\pi_g}
& \Gamma_{g,U}\backslash\Omega^r \ar@{^{ (}->}[r]^-{\theta_g} 
& \CU_g  \ar[d]^{\bar\pi_g} \\
M^r_{A,K}(\Cinf) \ar@{^{ (}->}[rr] && \OM^r_{A,K}(\Cinf). \\}$$
\item[(b)] This morphism is \'etale and its image is a Zariski open subset of $\OM^r_{A,K}(\Cinf)$. 
\item[(c)] Varying $g\in\GL_r(\BAfF)$, the union of the images of the different maps $\bar\pi_g$ is equal to $M^{r,+}_{A,K}(\Cinf)$ for a certain Zariski open subset $M^{r,+}_{A,K}$ of $\OM^r_{A,K}$ whose complement has codimension $\ge2$.
\end{itemize}
\end{Prop}

\begin{Proof}
This is due to Kapranov \cite{Kapranov} in the special case $A=\BF_q[t]$, and to H\"aberli \cite{Haeberli} in the general case. 
\end{Proof}

\begin{Rem}\label{AvoidKapranovHaeberli}
\rm For our application of Proposition \ref{Codim1Prop1} in the proof of Lemma \ref{AAMFLem2}, it would suffice to have, for every~$g$, an \'etale morphism on some arbitrarily small open subset $\CV_g\subset\CU_g$ that is not contained in $\Cinf^\times\times\Omega^{r-1}$, such that every connected component of codimension $1$ of $\OM^r_{A,K}(\Cinf) \setminus M^r_{A,K}(\Cinf)$ contains a point in the image of $\CV_g$ for some~$g$. It is probably possible to prove this without the explicit description of $\smash{\OM^r_{A,K}(\Cinf)}$ by Kapranov and H\"aberli, using only the fact from \cite[\PropositionCite 4.10]{Pink} that the fiber of the universal family $(\bar{E},\bar{\phi})$ over the generic point of any irreducible component of codimension $1$ of $\OM^r_{A,K} \setminus M^r_{A,K}$ is a Drinfeld $A$-module of rank $r-1$. But it would be a shame not to use the wonderful results from \cite{Kapranov} and \cite{Haeberli} when they are available.
%
%
\end{Rem}

Next let $(\BGacomma{\CU_g},\tilde\psi^{L_g})$ be the generalised Drinfeld $A$-module over $\CU_g$ that is furnished by Proposition \ref{GenDrinModCodim1}.

\begin{Prop}\label{Codim1Prop2}
There exists a unique isomorphism of generalised Drinfeld modules over~$\CU_g$
$$\bar\pi_g^*(\bar{E},\bar{\phi})\ \stackrel{\sim}{\longto}\ 
(\BGacomma{\CU_g},\bar\psi^{L_g}),$$
whose pullback under $\theta_g\circ\pi_{\Gamma_{g,U}}:$ $\Omega^r\to\CU_g$ is the isomorphism 
$$\xymatrix@C+10pt{
\pi_{\Gamma_{g,U}}^*\theta_g^*\bar\pi_g^*(\bar{E},\bar{\phi})
\ar@{=}[r]^-{\rm \ref{Codim1Prop1}\;(a)} & \pi_g^*(E,\phi)
\ar[r]_-\sim^-{(\ref{Tuvok})} & (\BGacomma{\Omega^r},\psi^{L_g})
\ar[r]_-\sim^-{\ref{GenDrinModCodim1}} & \pi_{\Gamma_{g,U}}^*\theta_g^*(\BGacomma{\CU_g},\bar\psi^{L_g}).\\}$$
\end{Prop}

\begin{Proof}
Over $\CU_g\cap(\Cinf^\times\times\Omega^{r-1})$ the isomorphism is obtained from the construction preceding (\ref{PiLDrinU}). The extension to $\CU_g$ follows from analytic versions of \cite[\PropositionsCite 3.7--8]{Pink}, which say that homomorphisms and isomorphisms of generalised Drinfeld modules extend uniquely under open dense embeddings of normal integral schemes, and whose proofs work equally well in the analytic setting.
\end{Proof}


\begin{Prop}\label{Func3}
In the situation of Proposition \ref{Func1} we have:
\begin{itemize}
\item[(a)] The morphism $J_h: M^r_{A,K'}\to M^r_{A,K}$ extends uniquely to a morphism $\bar J_h: {\OM^r_{A,K'}\to \OM^r_{A,K}}$.
\end{itemize}
\medskip\noindent Now assume that $K$ and $K'$ are fine, and let $(\bar E,\bar\phi)$ and $(\bar E',\bar\phi')$ denote the respective universal families on $\OM^r_{A,K}$ and $\OM^r_{A,K'}$. Then:
\begin{itemize}
\item[(b)] If $h$ has coefficients in~$\hatA$, the isogeny $\eta_h: J_h^*(E,\phi) \to (E',\phi')$ extends uniquely to an isogeny $\bar\eta_h: \bar J_h^*(\bar E,\bar\phi) \to (\bar E',\bar\phi')$.
\item[(c)] If $h^{-1}$ has coefficients in~$\hatA$, the isogeny $\xi_h: (E',\phi')\to J_h^*(E,\phi)$ extends uniquely to an isogeny $\bar\xi_h: (\bar E',\bar\phi')\to \bar J_h^*(\bar E,\bar\phi)$.
\end{itemize}
\end{Prop}

\begin{Proof}
(Sketch) Assertions (a) and (c) are proved in \cite[\PropositionCite4.11]{Pink}. The same kinds of arguments establish (b).
\end{Proof}

\medskip
Finally, the formulas in Proposition \ref{Func1} (e), (f) and in Proposition \ref{Func2} automatically extend to the respective Satake compactification, because the extended morphisms already exist and two morphisms on an integral scheme are equal if they coincide on an open dense subscheme.

\section{Analytic versus algebraic modular forms}
\label{Sec:AAMF}

We keep the notation from the preceding section, and first we also assume that $K$ is fine. 
Let $\Lie\bar{E}$ denote the Lie algebra of~$\bar{E}$, which is an invertible coherent sheaf of modules on $\smash{\OM^r_{A,K}}$. (It is naturally isomorphic to the sheaf of sections of~$\bar{E}$, but in the present context it is safer to view it as the Lie algebra.) Consider the dual invertible sheaf $\CL := (\Lie\bar{E})^\vee$. By \cite[\TheoremCite 5.3]{Pink} this is ample. For any integer $k$ we abbreviate $\CL^k := \CL^{\otimes k}$. Following \cite[\DefinitionCite 5.4]{Pink} we have:

\begin{Def}\label{AlgModFormsDef}
An \emph{algebraic Drinfeld modular form of weight $k$ and level $K$} is an element of the space
$$\Malg_k(M^r_{A,K})\ :=\ H^0(\OM^r_{A,K},\CL^k).$$
\end{Def}

Since $\OM^r_{A,K}$ is a projective algebraic variety with field of constants 
$F_K$, this is a finite-dimensional vector space over~$F_K$ or, depending on one's point of view, over~$F$. Our aim is to relate it with a space of analytic modular forms. Note that the decomposition (\ref{AnalUniv4a}) yields natural isomorphisms
\UseTheoremCounterForNextEquation
\begin{equation}\label{ModMapAll}
\Malg_k(M^r_{A,K})\otimes_F\Cinf\ \cong\ 
H^0(\OM^r_{A,K}\times_{\Spec F}\Spec\Cinf,\CL^k)
\ \cong\ \smash{\bigoplus_{i=1}^n H^0(\OM_{g_i},\CL^k).}
\end{equation}
Also, any irreducible component $\OM_g$ of \smash{$\OM^r_{A,K}\times_{\Spec F}\Spec\Cinf$} has field of definition $F_K$; hence pullback induces an isomorphism 
\UseTheoremCounterForNextEquation
\begin{equation}\label{ModMap1}
\Malg_k(M^r_{A,K}) \otimes_{F_K}\Cinf
\ \cong\ H^0(\OM_g,\CL^k).
\end{equation}
Let $\CL^\an$ denote the invertible sheaf on the rigid analytic space $\OM^r_{A,K}(\Cinf)$ obtained from~$\CL$. Its pullback $\pi_g^*\CL^\an$ is an invertible sheaf on~$\Omega^r$, which must be trivial, because $\Omega^r$ is a Stein space (\cite[\PropositionCite 4]{SchneiderStuhler}).
In fact, we have an explicit trivialisation: The isomorphism of line bundles $\pi_g^* E \to \BGacomma{\Omega^r}$ underlying the isomorphism of Drinfeld modules (\ref{Tuvok}) induces an isomorphism for the dual of the sheaf of sections 
\UseTheoremCounterForNextEquation
\begin{equation}\label{InvSheafIdent}
\pi_g^* \CL^\an\ \stackrel{\sim}{\longto}\ \CO_{\Omega^r}.
\end{equation}
Via this trivialisation, the pullback of any section $s\in H^0(M^r_{A,K}(\Cinf),(\CL^\an)^k)$ becomes a holomorphic function $\pi_g^* s:\Omega^r\to\Cinf$.

\begin{Lem}\label{AAMFLem0}
For any section $s\in H^0(M^r_{A,K}(\Cinf),(\CL^\an)^k)$ and any $g\in\GL_r(\BAfF)$ and $\gamma\in\GL_r(F)$ and $k\in K$ we have 
$$\pi_g^* s\ =\ (\pi_{\gamma gk}^*s)|_k\gamma.$$
\end{Lem}

\begin{Proof}
Since $\CL$ is the dual of the invertible sheaf of sections of $\bar{E}$, the commutative diagram (\ref{AnalUnivPullback2a}) yields a commutative diagram
$$\qquad\xymatrix@C-10pt{
\ \pi_g^*(\CL^\an)^k\ 
\ar[rrrrr]_-\sim^-{(\ref{InvSheafIdent})\;{\rm for}\;g}
\ar@{=}[d]^{(\ref{AnalUnivPig2})}
&&&&& \ \CO_{\Omega^r}\ 
\ar[d]_\wr^{{\rm multiplication\ by\ } j(\gamma,\underline{\ })^k}\\
\ \gamma^* \pi_{\gamma gk}^*(\CL^\an)^k\ 
\ar[rrrr]_-\sim^-{(\ref{InvSheafIdent})\;{\rm for}\;\gamma gk}
&&&& \ \gamma^*\CO_{\Omega^r}\ \ar@{=}[r]
& \ \CO_{\Omega^r}\rlap{.}\ \\}$$
For any $\omega\in\Omega^r$, evaluating $s$ at the point $\pi_g(\omega)=\pi_{\gamma gk}(\gamma(\omega))$ 
therefore yields the equality
$$j(\gamma,\omega)^k \cdot (\pi_g^*s)(\omega)\ =\ 
(\pi_{\gamma gk}^*s)(\gamma(\omega)).$$
In view of (\ref{DefOfActionOnF}) this implies that
$$(\pi_g^*s)(\omega)\ =\ 
j(\gamma,\omega)^{-k} \cdot (\pi_{\gamma gk}^*s)(\gamma(\omega))\ =\ 
((\pi_{\gamma gk}^*s)|_k\gamma)(\omega),$$
as desired.
\end{Proof}

\begin{Lem}\label{AAMFLem1}
The map $\pi_{g}^*$ induces an isomorphism
$$H^0(M_g(\Cinf),(\CL^\an)^k)\ \stackrel{\sim}{\longto}\ \CW_k(\Gamma_g).$$
\end{Lem}

\begin{Proof}
By definition the pullback by $\pi_g$ yields an isomorphism from $H^0(M_g(\Cinf),(\CL^\an)^k)$ to the space of $\Gamma_g$-invariant sections in $H^0(\Omega^r,\pi_g^*(\CL^\an)^k)$. 
But for every $\gamma\in\Gamma_g$ we have $\pi_{\gamma g}=\pi_g\circ\gamma^{-1} = \pi_g$ by (\ref{AnalUnivPig2}); so by Lemma \ref{AAMFLem0} the $\gamma$-invariance translates into the formula $\pi_g^*s = (\pi_g^*s)|_k\gamma$. By Definition \ref{Def:WeakModForm} the image of $\pi_g^*$ is therefore just the space of weak modular forms $\CW_k(\Gamma_g)$.
\end{Proof}

\begin{Lem}\label{AAMFLem2}
The map $\pi_g^*$ induces an isomorphism
$$H^0(\OM_g,\CL^k)\ \stackrel{\sim}{\longto}\ \CM_k(\Gamma_g).$$
\end{Lem}

\begin{Proof}
By rigid analytic GAGA due to K\"opf \cite[Satz\;4.7]{Koepf}, analytification yields an isomorphism $H^0(\OM_g,\CL^k) \stackrel{\sim}{\to} H^0(\OM_g(\Cinf),(\CL^\an)^k)$. 
Next, set $M_g^+ := \OM_g\cap M^{r,+}_{A,K}(\Cinf)$ for the Zariski open subset $M^{r,+}_{A,K}$ of $\OM^r_{A,K}$ from Proposition \ref{Codim1Prop1} (c). Since $\OM_g$ is normal integral and the complement $\OM_g\setminus M_g^+$ has codimension $\ge2$, by Bartenwerfer \cite[Satz\;10]{Bartenwerfer} the restriction map induces an isomorphism $H^0(\OM_g(\Cinf),(\CL^\an)^k) \stackrel{\sim}{\to} H^0(M_g^+(\Cinf),(\CL^\an)^k)$. 
By Lemma \ref{AAMFLem1} any section $s\in H^0(M_g(\Cinf),(\CL^\an)^k)$ corresponds to a weak modular form $\pi_g^*s\in\CW_k(\Gamma_g)$. It remains to determine when $s$ extends to a section in $H^0(M_g^+(\Cinf),(\CL^\an)^k)$.

We first analyse when it extends to the image of the map $\bar\pi_g$ from  Proposition \ref{Codim1Prop1}~(a). 
Recall that $\CL$ was defined as the dual of the invertible sheaf of sections of~$\bar{E}$. Thus the isomorphism of generalised Drinfeld modules in Proposition \ref{Codim1Prop2} induces an isomorphism
\UseTheoremCounterForNextEquation
\begin{equation}\label{InvSheafIdent2}
\bar\pi_g^*\CL^\an\ \cong\ \CO_{\CU_g}.
\end{equation}
Let $\bar\theta: \Omega^r\to\CU_g$ be the composite morphism in the top row of the diagram in Proposition \ref{Codim1Prop1} (a). Then by construction the pullback of the trivialisation (\ref{InvSheafIdent2}) to $\Omega^r$ via $\bar\theta$ is just the trivialisation in (\ref{InvSheafIdent}). 
Thus $s$ extends to a section of $(\CL^\an)^k$ over the image of $\bar\pi_g$ if and only if the function $\pi_g^*s:\Omega^r\to\Cinf$ is the pullback via $\bar\theta$ of a holomorphic function $\CU_g\to\Cinf$.
Here $\pi_g^*s$ is already a $\Gamma_U$-invariant function and therefore possesses a $u$-expansion by Proposition \ref{ModFormsLaurentExpansion}. Thus it is the pullback of a holomorphic function on~$\CU_g$ if and only if it is holomorphic at infinity in the sense of Definition \ref{Def:OrderAtInfinity}.

Now recall that for any $g,g'\in\GL_r(\BAfF)$ we have $M_g=M_{g'}$ if and only if $g'=\gamma gk$ for some $\gamma\in\GL_r(F)$ and $k\in K$. By Proposition \ref{Codim1Prop1} (c) the partial compactification $M_g^+$ is therefore the union of the images of the maps $\bar\pi_{\gamma gk}$ for all such $\gamma$ and~$k$. By the above argument for $\gamma gk$ in place of~$g$, it follows that $s$ extends to a section in $H^0(M_g^+(\Cinf),(\CL^\an)^k)$ if and only if for all $\gamma$ and $k$ the pullback $\pi_{\gamma gk}^*s$ is holomorphic at infinity. But by Lemma \ref{AAMFLem0} we have $\pi_{\gamma gk}^*s = (\pi_g^*s)|_k\gamma^{-1}$. Varying $\gamma$ we thus conclude that $\pi_g^*$ induces an isomorphism from $H^0(M_g^+(\Cinf),(\CL^\an)^k)$ to the space of modular forms $\CM_k(\Gamma_g)$. Combining everything yields the desired result.
\end{Proof}

\begin{Thm}\label{AAMFThm1}
If $K$ is fine, the maps $\pi_g^*$ and the isomorphisms (\ref{ModMap1}) respectively (\ref{ModMapAll}) induce isomorphisms
$$\Malg_k(M^r_{A,K})\otimes_{F_K}\Cinf\ \stackrel{\sim}{\longto}\ 
\CM_k(\Gamma_g),$$
$$\Malg_k(M^r_{A,K})\otimes_F\Cinf\ \ \stackrel{\sim}{\longto}\ \ 
\smash{\bigoplus_{i=1}^n \CM_k(\Gamma_{g_i}).}$$
\end{Thm}

\begin{Proof}
Direct consequence of Lemma \ref{AAMFLem2}.
\end{Proof}

\medskip
The above isomorphisms are functorial in the following sense. Consider a second fine open compact subgroup $K'<\GL_r(\hatA)$ and an element $h\in\GL_r(\BAfF)$ such that $hK'h^{-1}<K$. By Proposition \ref{Func3} (a) this data determines a morphism $\bar J_h: \OM^r_{A,K'}\longto \OM^r_{A,K}$. As before let $(\bar E',\bar\phi')$ denote the universal generalised Drinfeld module on $\OM^r_{A,K'}$. Let $\CL'$ denote the dual of the invertible sheaf of sections of~$\bar E'$. 

With $h$ fixed, consider any sufficiently divisible scalar $a\in A\setminus\{0\}$, so that the element $ha\in\GL_r(\BAfF)$ has coefficients in~$\hatA$. As a consequence of Propositions \ref{Func1} (f) and \ref{Func2}, we then have $\bar J_{ha}=\bar J_h$. The derivative of the isogeny $\bar\eta_{ha}$ in Proposition \ref{Func3} (b) thus induces an isomorphism 
$$(d\bar\eta_{ha})^\vee:\ \bar J_h^*\CL = \bar J_{ha}^*\CL \stackrel{\sim}{\longto} \CL'.$$

\begin{Lem}\label{LFunc1}
The isomorphism
$$\rho_h := a\cdot(d\bar\eta_{ha})^\vee:\ \bar J_h^*\CL \stackrel{\sim}{\longto} \CL'$$
is independent of the choice of~$a$.
\end{Lem}

\begin{Proof}
Consider a second element $b\in A\setminus\{0\}$ such that $hb$ has coefficients in~$\hatA$. Then so does $hab$, and Propositions \ref{Func2} (b) and \ref{Func1} (f) imply that $\eta_{hab} = \eta_{b}\circ\eta_{ha} = \phi'_b\circ\eta_{ha}$. Taking derivatives we deduce that $d\eta_{hab} = d\phi'_b\circ d\eta_{ha} = b\cdot d\eta_{ha}$ and hence $ab\cdot(d\eta_{hab})^\vee = ab\cdot b^{-1}\cdot (d\eta_{ha})^\vee = a\cdot (d\eta_{ha})^\vee$. Interchanging $a$ and $b$ implies that $ab\cdot(d\eta_{hab})^\vee = b\cdot (d\eta_{hb})^\vee$ and  hence $a\cdot (d\eta_{ha})^\vee = b\cdot (d\eta_{hb})^\vee$. Finally, this equality over the dense open subscheme $M^r_{A,K'}$ automatically extends to an equality over $\OM^r_{A,K'}$.
\end{Proof}

\medskip
Using pullback and the isomorphism $\rho_h$ we can now define a natural $F$-linear pullback map on modular forms, again denoted~$J_h^*$, by the commutative diagram
\UseTheoremCounterForNextEquation
\begin{equation}\label{LFunc2}
\hskip20pt\vcenter{\xymatrix@R-15pt@C+5pt{
\llap{$J_h^*:\ \ $} \Malg_k(M^r_{A,K}) \ar[rr] && \Malg_k(M^r_{A,K'}) \\
H^0(\OM^r_{A,K},\CL^k) \ar@{=}[u] \ar[r]^-{\bar J_h^*} & 
H^0(\OM^r_{A,K'},\bar J_h^*\CL^k) \ar[r]^-{\rho_h^k} & 
H^0(\OM^r_{A,K'},\CL^{\prime k})\rlap{.} \ar@{=}[u] \\}}
\end{equation}
To describe its behavior under the isomorphisms from Theorem \ref{AAMFThm1}, for any $g\in\GL_r(\BAfF)$ consider the arithmetic subgroup $\Gamma'_{gh} := \GL_r(F)\cap ghK'(gh)^{-1}$, which by construction is contained in the arithmetic subgroup $\Gamma_g := \GL_r(F)\cap gKg^{-1}$.

\begin{Prop}\label{LFunc3}
For any $g\in\GL_r(\BAfF)$ the diagram
$$\xymatrix@C+10pt{
\Malg_k(M^r_{A,K}) \ar[r]^-{J_h^*} \ar[d]^{\pi_g^*}
& \Malg_k(M^r_{A,K'}) \ar[d]^{\pi_{gh}^{\prime*}} \\
\CM_k(\Gamma_g) \ar@{^{ (}->}[r] & \CM_k(\Gamma'_{gh}) \\}$$
commutes, where the horizontal map on the bottom is the inclusion map.
\end{Prop}

\begin{Proof}
Assume first that $h$ has coefficients in~$\hatA$. As $\CL$ and $\CL'$ are the duals of the invertible sheaves of sections of $\bar{E}$ and $\bar{E}'$, Proposition \ref{Func1} (c) yields a commutative diagram
$$\xymatrix@C-16pt{
\pi_g^*\CL^\an \ar@{=}[r] \ar@{=}[d]_\wr^{(\ref{InvSheafIdent})\;{\rm for}\;g} &
\pi_{gh}^{\prime*}J_h^*\CL^\an \ar[rrrrrrrr]^{\pi_{gh}^{\prime*}\rho_h\;=\;\pi_{gh}^{\prime*}(d\eta_h)^\vee} &&&&&&&&
\pi_{gh}^{\prime*}\CL^{\prime\an} \ar@{=}[d]_\wr^{(\ref{InvSheafIdent})\;{\rm for}\;gh} \\
\CO_{\Omega^r} \ar[rrrrrrrrr]^{(d\tilde\eta_h)^\vee} &&&&&&&&& 
\CO_{\Omega^r}.\\}$$
By the construction (\ref{xih1}) of $\tilde\eta_h$ we have $d\tilde\eta_h=1$. 
The desired commutativity thus follows from the definition of $\pi_g^*$ and $\pi_{gh}^{\prime*}$.

In the general case take any $a\in A\setminus\{0\}$ such that $ha\in\GL_r(\BAfF)$ has coefficients in~$\hatA$. Repeating the above calculation twice with $(g,h)$ replaced by $(g,ha)$ and $(gh,a)$, respectively, and noting that $\pi'_{gha}=\pi'_{gh}$, yields a commutative diagram
$$\xymatrix@C-16pt{
\pi_g^*\CL^\an \ar@{=}[d]_\wr^{(\ref{InvSheafIdent})\;{\rm for}\;g}
 \ar[rrrrrr]^{\pi_{gha}^{\prime*}(d\eta_{ha})^\vee} &&&&&&
\pi_{gha}^{\prime*}\CL^{\prime\an} \ar@{=}[d]_\wr^{(\ref{InvSheafIdent})\;{\rm for}\;gha}  
&&&&&&
\ar[llllll]_{\pi_{gha}^{\prime*}(d\eta_{a})^\vee} 
\pi_{gh}^{\prime*}\CL^{\prime\an} \ar@{=}[d]_\wr^{(\ref{InvSheafIdent})\;{\rm for}\;gh}  \\
\CO_{\Omega^r} \ar[rrrrrr]^{\id} &&&&&&
\CO_{\Omega^r} &&&&&& \CO_{\Omega^r}. \ar[llllll]_{\id}\\}$$
Here $d\eta_{a}=d\phi'_a=a$ by Proposition \ref{Func1} (f), hence the upper horizontal arrow on the right is multiplication by $a^{-1}$. Together we thus obtain the commutative diagram
$$\xymatrix@C-10pt{
\pi_g^*\CL^\an \ar@{=}[d]_\wr^{(\ref{InvSheafIdent})\;{\rm for}\;g}
 \ar[rrrrrr]^{\pi_{gh}^{\prime*}\rho_h \;=\; a\cdot\pi_{gha}^{\prime*}(d\eta_{ha})^\vee} &&&&&&  
\pi_{gh}^{\prime*}\CL^{\prime\an} \ar@{=}[d]_\wr^{(\ref{InvSheafIdent})\;{\rm for}\;gh}  \\
\CO_{\Omega^r} \ar[rrrrrr]^{\id} &&&&&& \CO_{\Omega^r}, \\}$$
and again the desired commutativity follows from the definition of $\pi_g^*$ and $\pi_{gh}^{\prime*}$.
\end{Proof}

\begin{Prop}\label{LFunc4}
\begin{itemize}
\item[(a)] If $K=K'$ and $h\in K$, then $J_h^*=\id$.
\item[(b)] If $K=K'$ and $h=a\cdot\Id_r$ for $a\in A\setminus\{0\}$ then $J_h^*=a^k\cdot\id$.
\item[(c)] For any fine open compact subgroups $K,K',K''<\GL_r(\hatA)$ and elements $h,h'\in\GL_r(\BAfF)$ such that $hK'h^{-1}<K$ and $h'K''h^{\prime-1}<K'$, we have $J_{hh'}^*=J_{h'}^*\circ J_h^*$.
\end{itemize}
\end{Prop}

\begin{Proof}
Direct computation using Proposition \ref{Func2}.
\end{Proof}

\medskip
Now recall that the elements $g_1,\ldots,g_n$ appearing in Theorem \ref{AAMFThm1} are the representatives of the double quotient $\GL_r(F)\backslash\GL_r(\BAfF)/K$ used in (\ref{AnalUnivDec}). Likewise choose representatives $g'_1,\ldots,g'_{n'}$ of the double quotient $\GL_r(F)\backslash\GL_r(\BAfF)/K'$. For each $1\le j\le n'$ consider the arithmetic subgroup $\Gamma'_{g'_j} := \GL_r(F)\cap g'_jK'g_j^{\prime-1}$, and choose $1\le i_j\le n$ and $\gamma_j\in\GL_r(F)$ and $k_j\in K$ such that $\gamma_jg'_jh^{-1}k_j=g_{i_j}$. Then direct calculations show that $\gamma_j\Gamma'_{g'_j}\gamma_j^{-1} < \Gamma_{g_{i_j}}$ and that the following diagram commutes:
\vskip-10pt
\UseTheoremCounterForNextEquation
\begin{equation}\label{AnalUnivFunct}
\vcenter{\xymatrix@C+5pt{
\smash{\displaystyle\coprod_{j=1}^{\;n'}}\, \Gamma'_{g'_j}\backslash\Omega^r \ar[r]_-\sim^-{(\pi'_{\smash{\scriptscriptstyle g'_j}})} \ar[d]
& \GL_r(F)\!\bigm\backslash\! \bigl( \Omega^r\times\GL_r(\BAfF)/K'\bigr) \ar[r]^-\sim \ar[d]
& M^r_{A,K'}(\Cinf) \ar[d]^{J_h} \\
\smash{\displaystyle\coprod_{i=1}^n}\, \Gamma_{g_i}\backslash\Omega^r \ar[r]_-\sim^-{(\pi_{g_i})}
& \GL_r(F)\!\bigm\backslash\! \bigl( \Omega^r\times\GL_r(\BAfF)/K\bigr) \ar[r]^-\sim
& M^r_{A,K}(\Cinf)\rlap{,} \\}}
\end{equation}
where the vertical map in the middle is $[(\omega,g)]\mapsto[(\omega,gh^{-1})]$ and the one on the left sends a coset $\Gamma'_{g'_j}\omega$ in the $j$-th subset to the coset $\Gamma_{g_{i_j}}\gamma_j(\omega)$ in the $i_j$-th subset.

\begin{Prop}\label{AnalUnivFunctLem}
If $K$ and $K'$ are fine, the map $J_h^*$ from (\ref{LFunc2}) and the isomorphisms from Theorem \ref{AAMFThm1} for $K'$ and $K$ fit into a commutative diagram
$$\xymatrix@C+15pt@R-15pt{
\Malg_k(M^r_{A,K})\otimes_F\Cinf 
\ar[dd]_\wr^{\ref{AAMFThm1}}
\ar[r]^-{J_h^*\otimes\id} 
& \Malg_k(M^r_{A,K'})\otimes_F\Cinf 
\ar[dd]_\wr^{\ref{AAMFThm1}} \\ &&\\
\smash{\;\displaystyle\bigoplus_{i=1}^n\;} \CM_k(\Gamma_{g_i}) \ar[r] 
& \smash{\;\displaystyle\bigoplus_{j=1}^{\;n'}\;} \CM_k(\Gamma'_{g'_j}) \\
(f_i)_{i=1}^n \ar@{|->}[r] & 
(f_{i_j}|_k\gamma_j)_{j=1}^{n'}\rlap{.} \\}$$
\end{Prop}

\begin{Proof}
For each $1\le j\le n'$ we have a commutative diagram
$$\xymatrix@C-20pt{
& \Malg_k(M^r_{A,K}) \ar[rrrr]^-{J_h^*} \ar[dr]^>>>>{\pi_{g'_jh^{-1}}^*}
\ar[dl]_>>>>{\pi_{g_{i_j}}^*}
&&&& \Malg_k(M^r_{A,K'}) \ar[d]^{\pi_{g'_j}^{\prime*}} \\
\CM_k(\Gamma_{g_{i_j}}) \ar[rr]^{f\;\mapsto\; f|_k\gamma_j}
&& \CM_k(\Gamma_{g'_jh^{-1}}) \ar@{^{ (}->}[rrr]^{\rm incl.} 
&&& \CM_k(\Gamma'_{g'_j})\rlap{,} \\}$$
which commutes on the left by the equation $\gamma_jg'_jh^{-1}k_j=g_{i_j}$ and Lemma \ref{AAMFLem0}, and on the right by Proposition \ref{LFunc3} for $g=g'_jh^{-1}$. Summing over all $j$ yields the desired formula.
\end{Proof}

\medskip
Finally consider an arbitrary open compact subgroup $K<\GL_r(\BAfF)$. Let $\tilde K$ be any open normal subgroup of $K$ which is fine, for instance, the principal congruence subgroup $K(N)$ for a sufficiently divisible non-zero ideal $N\subsetneqq A$. Then by Proposition \ref{LFunc4} the maps $J_h^*$ for all $h\in K$ induce a right action of $K/\tilde K$ on the space of modular forms of level~$\tilde K$. In \cite[\DefinitionCite5.4]{Pink} we defined:

\begin{Def}\label{AlgModFormsDef2}
The space of \emph{algebraic Drinfeld modular forms of weight $k$ and arbitrary level $K$} is the space of $K$-invariants
$$\Malg_k(M^r_{A,K})\ :=\ \Malg_k(M^r_{A,\tilde K})^{K}.$$
\end{Def}

Once defined using one choice of~$\tilde K$, the same equality then holds for arbitrary open compact subgroups $\tilde K\triangleleft K<\GL_r(\BAfF)$. This makes $\Malg_k(M^r_{A,K})$ independent of the choice of~$\tilde K$.
Moreover, for any $g\in\GL_r(\BAfF)$ we define the pullback map $\pi_g^*$ on $\Malg_k(M^r_{A,K})$ as the restriction of the map $\pi_g^*$ on $\Malg_k(M^r_{A,\tilde K})$. Using Proposition \ref{LFunc3} in the case $h=\Id_r$ we find that this is again independent of the choice of~$\tilde K$. Likewise we can define a map $J_h^*: \Malg_k(M^r_{A,K}) \to \Malg_k(M^r_{A,K'})$ for arbitrary $h$, $K$, $K'$ as the restriction to $K$-, resp. $K'$-invariants from suitable smaller open compact subgroups. With this we can now conclude:

\begin{Prop}\label{AllHoldsForNotFineK}
Theorem \ref{AAMFThm1} and Propositions \ref{LFunc3} and \ref{LFunc4} and \ref{AnalUnivFunctLem} hold for arbitrary open subgroups.
\end{Prop}

\begin{Proof}
(Sketch) For all $h\in K$ we have $hK'h^{-1}=K'$, so using Proposition \ref{AnalUnivFunctLem} with $K$ replaced by $K'$ we can translate the right action of $K/K'$ on $\Malg_k(M^r_{A,K'})\otimes_F\Cinf$ to the space $\bigoplus_{j=1}^{\;n'} \CM_k(\Gamma'_{g'_j})$. This action interchanges the summands $\CM_k(\Gamma'_{g'_j})$ whenever $g'_j$ lies in the same coset $\GL_r(F)g_iK$, and the stabiliser of such a summand acts through the action of all $\gamma\in\Gamma_{g_i}$ by $f\mapsto f|_k\gamma$. But the space of invariants in $\CM_k(\Gamma'_{g'_j})$ under this action is simply $\CM_k(\Gamma_{g_i})$. Taking invariants we thus deduce the second isomorphism in Theorem \ref{AAMFThm1} for the group~$K$.
The remaining statements follow in the same way by taking invariants in each case.
\end{Proof}

\section{Finiteness results}
\label{Sec:Fin}

\begin{Thm}\label{FiniteDimension}
For any congruence subgroup $\Gamma<\GL_r(F)$ we have:
\begin{itemize}
\item[(a)] $\dim_{\Cinf}\CM_{k,m}(\Gamma)<\infty$ for any integers $k$ and~$m$.
\item[(b)] $\CM_{k,m}(\Gamma)=0$ whenever $k<0$ and $r\ge2$.
\item[(c)] The graded ring $\CM_*(\Gamma) := \bigoplus_{k\ge0}\CM_{k}(\Gamma)$ is a normal integral domain that is finitely generated as a $\Cinf$-algebra.
\end{itemize}
\end{Thm}

\begin{Proof}
First assume that $\Gamma$ is the principal congruence subgroup $\Gamma(N)$ associated to some level $0\not=N\subsetneqq A$. Setting $K:=K(N)$, for $g=1$ the arithmetic subgroup $\Gamma_g$ from (\ref{AnalUnivLion}) is then~$\Gamma$. By Theorem \ref{AAMFThm1} we thus have $H^0(\OM^r_{A,K},\CL^k)\otimes_{F_K}\Cinf \cong \CM_k(\Gamma)$. As space of sections of a coherent sheaf on a projective algebraic variety it is therefore finite dimensional, proving (a). Moreover, since $\CL$ is ample by \cite[\TheoremCite5.3]{Pink}, this space is zero if ${k<0}$ and every irreducible component of the variety has dimension $\ge1$, proving (b). Also, the ring $\bigoplus_{k\ge0}H^0(\OM^r_{A,K},\CL^k)$ is a normal integral domain that is finitely generated as a $\Cinf$-algebra by \cite[\TheoremCite5.6]{Pink}, proving (c). 

Next, for any two congruence subgroups $\Gamma'\triangleleft\Gamma$, the respective space or graded ring for $\Gamma$ is obtained from that for $\Gamma'$ by taking invariants under a certain action of the finite group $\Gamma/\Gamma'$. The statements for $\Gamma$ thus follow from those for~$\Gamma'$.

Finally, for an arbitrary congruence subgroup $\Gamma<\GL_r(F)$ consider the finitely generated $A$-submodule $L := \Gamma\cdot A^r\subset F^r$, and choose an ideal $0\not=I\subsetneqq A$ such that $IL\subset A^r$. Let $\Gamma'$ be the subgroup of elements of $\Gamma$ that act trivially on $L/IL$. Then $\Gamma'\triangleleft\Gamma$ and $\Gamma'<\GL_r(A)$. Also $\Gamma'$ is again a congruence subgroup, so it contains $\Gamma(N)$ for some level $0\not=N\subsetneqq A$. As $\Gamma'<\GL_r(A)$, we then have $\Gamma(N)\triangleleft\Gamma'\triangleleft\Gamma$, and the statements for $\Gamma$ follow from those for $\Gamma(N)$ by applying the above reduction step twice.
\end{Proof}

\begin{Prop}\label{CuspEx}
Let $\Gamma<\GL_r(A)$ be a congruence subgroup whose image in $\GL_r(A/\Fp)$ is unipotent for some maximal ideal $\Fp\subset A$. Then for every $k\gg0$ there exists a non-zero cusp form of weight $k$ for~$\Gamma$.
\end{Prop}

For an explicit construction of such cusp forms using Eisenstein series see Remark \ref{NonZeroCuspForm}.

\medskip
\begin{Proof}
Choose a level $0\not=N\subsetneqq A$ such that $\Gamma(N)<\Gamma$, and set $K:=K(N)\cdot\Gamma<\GL_r(\hatA)$. Then $K$ is fine, and for $g=1$ we have $\Gamma_g=K\cap\GL_r(A)=\Gamma$. Let $\infty$ denote the reduced divisor on $\OM^r_{A,K}$ with support $\OM^r_{A,K}\smallsetminus M^r_{A,K}$. By Theorem \ref{AAMFThm1} and the definition of cusp forms we then have 
$$H^0(\OM^r_{A,K},\CL^k(-\infty))\otimes_{F_K}\Cinf 
\ \cong\ \CS_k(\Gamma).$$
As $\CL$ is ample, the left hand side is non-zero for all $k\gg0$, as desired.
\end{Proof}

\section{Hecke operators}
\label{HeckeOps}

Consider any element $h\in\GL_r(\BAfF)$ and any open compact subgroups $K$, $K'<\GL_r(\hat A)$ such that $hK'h^{-1}<K$. Then by (\ref{LFunc2}) and Proposition \ref{AllHoldsForNotFineK}, there is a well-defined \emph{pullback map}
\UseTheoremCounterForNextEquation
\begin{equation}\label{HO:Pull}
J_h^*\!:\ \Malg_k(M^r_{A,K})\ \longto\ \Malg_k(M^r_{A,K'})
\end{equation}
satisfying Proposition \ref{LFunc4}. 

We can also construct a natural map in the other direction. Since $J_h^*$ is an isomorphism if $hK'h^{-1}=K$, we restrict ourselves to the case that $h=\Id_r$ and $K'<K$. Choose an open subgroup $\tilde K<K'$ which is normal in~$K$. Then by Definition \ref{AlgModFormsDef2} we have
\UseTheoremCounterForNextEquation
\begin{equation}\label{HO:PullDiag}
\quad\vcenter{\xymatrix@R-10pt@C+20pt{
\Malg_k(M^r_{A,K})  \ar@{^{ (}->}[r]^-{J_{\Id_r}^*} \ar@{=}[d] & 
\Malg_k(M^r_{A,K'}) \ar@{^{ (}->}[dr]^-{J_{\Id_r}^*} \ar@{=}[d] & \\
\Malg_k(M^r_{A,\tilde K})^K \ar@{^{ (}->}[r] &
\Malg_k(M^r_{A,\tilde K})^{K'} \ar@{..>}@/^2ex/[l]^{\rm trace} 
\ar@{^{ (}->}[r] & \Malg_k(M^r_{A,\tilde K})  \rlap{.}\\}}
\end{equation}
We define the dotted arrow by
\UseTheoremCounterForNextEquation
\begin{equation}\label{HO:Trace}
f\longmapsto \mathop{\rm trace}(f) := \sum\nolimits_{h'} J_{h'}^*f,
\end{equation}
where $h'$ runs through a set of representatives of the quotient $K'\backslash K$. The composite of this trace map with the vertical isomorphisms in 
(\ref{HO:PullDiag}) is the \emph{pushforward map}
\UseTheoremCounterForNextEquation
\begin{equation}\label{MF:Push}
J_{\Id_r,*}\!:\ \Malg_k(M^r_{A,K'})\ \longto\ \Malg_k(M^r_{A,K}).
\end{equation}

\medskip
Now consider any element $h\in\GL_r(\BAfF)$ and any open compact subgroup $K<\GL_r(\hat A)$, bearing no particular relation with each other. Then 
we call the pair of morphisms
\UseTheoremCounterForNextEquation
\begin{equation}\label{HO:HC}
\xymatrix{
M^r_{A,K} & \ar[l]_-{J_h} M^r_{A,K \cap\nobreak h^{-1}Kh} \ar[r]^-{J_{\Id_r}} & M^r_{A,K} \\}
\end{equation}
the \emph{Hecke correspondence on $M^r_{A,K}$ associated to~$h$}. The composite map
\UseTheoremCounterForNextEquation
\begin{equation}\label{HO:Hecke}
\xymatrix{
\llap{$T_h\!:\ \ $}\Malg_k(M^r_{A,K}) \ar[r]^-{J^*_h} & 
\Malg_k(M^r_{A,K \cap\nobreak h^{-1}Kh}) \ar[r]^-{J_{\Id_r,*}} & \Malg_k(M^r_{A,K}) \\}
\end{equation}
is called the \emph{Hecke operator on $\Malg_k(M^r_{A,K})$ associated to~$h$}. It depends only on the double coset $KhK$.

\medskip
The composites of Hecke operators are calculated as follows:

\begin{Prop}
\label{HO:Comp}
For any $h$, $h'\in\GL_r(\BAfF)$ and any open compact subgroup $K<\GL_r(\hat A)$ the Hecke operators on $\Malg_k(M^r_{A,K})$ satisfy
$$T_{h'}\circ T_h \ =\ \sum_{h''} \; 
\bigl[ K\cap h^{\prime\prime-1}Kh'' : K\cap h^{-1}Kh\cap h^{\prime\prime-1}Kh'' \bigr]
\cdot T_{h''}$$
where $h''$ runs through a set of representatives of the double quotient
$$(hKh^{-1}\cap K) \backslash hKh' / (K\cap h^{\prime-1}Kh').$$
\end{Prop}

\begin{Proof}
This is \cite[Prop.\;6.10]{Pink} with the change of conventions taken into account.
\end{Proof}


\medskip
In the rest of this section we work out how the maps $J_{\Id_r,*}$ and $T_h$ translate under the isomorphism from Theorem \ref{AAMFThm1}.

\begin{Prop}\label{PushForwardDescription}
Consider any open compact subgroups $K'<K<\GL_r(\hat A)$ and any representatives $g_1,\ldots,g_n$ of the double quotient $\GL_r(F)\backslash\GL_r(\BAfF)/K$ and representatives $g'_1,\ldots,g'_{n'}$ of the double quotient $\GL_r(F)\backslash\GL_r(\BAfF)/K'$. For each $1\le i\le n$ consider the arithmetic subgroup $\Gamma_{g_i} := \GL_r(F)\cap g_iKg_i^{-1}$ and for each $1\le j\le n'$ the arithmetic subgroup $\Gamma'_{g'_j} := \GL_r(F)\cap g'_jK'g_j^{\prime-1}$. Then the map $J_{\Id_r,*}$ from (\ref{MF:Push}) and the isomorphisms from Theorem \ref{AAMFThm1} for $K'$ and $K$ fit into a commutative diagram
$$\qquad\xymatrix@C+20pt@R-15pt{
\Malg_k(M^r_{A,K'})\otimes_F\Cinf 
\ar[dd]_\wr^{\ref{AAMFThm1}}
\ar[r]^-{J_{\Id_r,*}\otimes\id} 
& \Malg_k(M^r_{A,K})\otimes_F\Cinf 
\ar[dd]_\wr^{\ref{AAMFThm1}} \\ &&\\
\smash{\;\displaystyle\bigoplus_{j=1}^{\;n'}\;} \CM_k(\Gamma'_{g'_j}) \ar[r] 
& \smash{\;\displaystyle\bigoplus_{i=1}^n\;} \CM_k(\Gamma_{g_i}) \rlap{,} \\
(f_j)_{j=1}^{n'} \ar@{|->}[r] & 
\smash{\;\Bigl(\displaystyle\sum_{j,\gamma}f_j|_k\gamma \Bigr){\strut}_{i=1}^n} \rlap{,} \\}$$
where, for each index~$i$, the sum extends over all pairs of indices $1\le j\le n'$ and elements $\gamma\in \GL_r(F)\cap g'_jKg_i^{-1}$ up to left multiplication by $\Gamma'_{g'_j}$.
\end{Prop}

\begin{Proof}
Suppose first that $K'\triangleleft K$. Then for any $h\in K$  and any $1\le i\le n$ there is an index $1\le j_{ih}\le n'$ and an element $\gamma_{ih}\in\GL_r(F)$ such that $g'_{j_{ih}}\in\gamma_{ih}g_ih^{-1}K'$. By Propositions \ref{AnalUnivFunctLem} and \ref{AllHoldsForNotFineK} the map $J_h^*\otimes\id$ thus corresponds to the map
$$(f_j)_{j=1}^{n'} \ \longmapsto\ (f_{j_{ih}}|_k\gamma_{ih} )_{i=1}^n.$$
Next observe that $j_{ih}$ is unique and $\gamma_{ih}$ is unique up to multiplication on the left by $\Gamma'_{g'_{j_{ih}}}$, and both depend only on $i$ and the coset $K'h$. Summing over all cosets $K'h\subset K$ thus shows that $J_{\Id_r,*}\otimes\id$ corresponds to the map
$$(f_j)_{j=1}^{n'} \ \longmapsto\ 
\sum_{K'h} (f_{j_{ih}}|_k\gamma_{ih} )_{i=1}^n 
\ =\ \Bigl(\displaystyle\sum_{j,\gamma}f_j|_k\gamma \Bigr)\strut_{i=1}^n$$
with the indicated summation over $(j,\gamma)$. This proves the assertion in the case $K'\triangleleft K$. 

In the general case, one must take an open compact subgroup $\tilde K<K'$ which is normal in~$K$, choose representatives for $\GL_r(F)\backslash\GL_r(\BAfF)/\tilde K$, write down the commutative diagrams from Proposition \ref{AnalUnivFunctLem} for the maps $J_{\Id_r}^*: \Malg_k(M^r_{A,K}) \to \Malg_k(M^r_{\smash{A,\tilde K}})$ and $J_{\Id_r}^*: \Malg_k(M^r_{A,K'}) \to \Malg_k(M^r_{\smash{A,\tilde K}})$ and $J_h^*\!:\ \Malg_k(M^r_{\smash{A,\tilde K}})\ \longto\ \Malg_k(M^r_{\smash{A,\tilde K}})$ for all $h\in K$, and eliminate everything concerning $\tilde K$ from the resulting expression for $J_{\Id_r,*}\otimes\id$. We leave this direct and tedious 
calculation to the reader.
\end{Proof}


\begin{Prop}\label{HeckeDescription}
Consider any element $h\in\GL_r(\BAfF)$, any open compact subgroup $K<\GL_r(\hatA)$ and any representatives $g_1,\ldots,g_n$ of the double quotient $\GL_r(F)\backslash\GL_r(\BAfF)/K$. Then the Hecke operator $T_h$ from (\ref{HO:Hecke}) and the isomorphism from Theorem \ref{AAMFThm1} fit into a commutative diagram
$$\qquad\xymatrix@C+20pt@R-15pt{
\Malg_k(M^r_{A,K})\otimes_F\Cinf 
\ar[dd]_\wr^{\ref{AAMFThm1}}
\ar[r]^-{T_h\otimes\id} 
& \Malg_k(M^r_{A,K})\otimes_F\Cinf 
\ar[dd]_\wr^{\ref{AAMFThm1}} \\ &&\\
\smash{\;\displaystyle\bigoplus_{i=1}^n\;} \CM_k(\Gamma_{g_i}) \ar[r] 
& \smash{\;\displaystyle\bigoplus_{i=1}^n\;} \CM_k(\Gamma_{g_i}) \rlap{,} \\
(f_i)_{i=1}^{n} \ar@{|->}[r] & 
\smash{\;\Bigl(\displaystyle\sum_{i',\delta}f_{i'}|_k\,\delta \Bigr){\strut}_{i=1}^n} \rlap{,} \\}$$
where, for each index~$i$, the sum extends over all pairs of indices $1\le i'\le n$ and elements $\delta\in \GL_r(F)\cap g_{i'}KhKg_i^{-1}$ up to left multiplication by $\Gamma_{g_{i'}}$. Moreover, the index $i'$ that actually occurs in the sum depends only $i$ and~$h$.
\end{Prop}

\begin{Proof}
Set $K':=K \cap\nobreak h^{-1}Kh$ and choose representatives $g'_1,\ldots,g'_{n'}$ of the double quotient $\GL_r(F)\backslash\GL_r(\BAfF)/K'$. For each $1\le j\le n'$ select an index $1\le i_j\le n$ and elements $\gamma_j\in \GL_r(F)$ and $k_j\in K$ such that $\gamma_jg'_jh^{-1}k_j=g_{i_j}$. Then by Propositions \ref{AnalUnivFunctLem} and \ref{PushForwardDescription} we have a commutative diagram
$$\qquad\xymatrix@C+20pt@R-15pt{
\Malg_k(M^r_{A,K})\otimes_F\Cinf 
\ar[dd]_\wr^{\ref{AAMFThm1}}
\ar[r]^-{J_h^*\otimes\id} 
& \Malg_k(M^r_{A,K'})\otimes_F\Cinf 
\ar[dd]_\wr^{\ref{AAMFThm1}}
\ar[r]^-{J_{\Id_r,*}\otimes\id} 
& \Malg_k(M^r_{A,K})\otimes_F\Cinf 
\ar[dd]_\wr^{\ref{AAMFThm1}} \\ &&&\\
\smash{\;\displaystyle\bigoplus_{i=1}^n\;} \CM_k(\Gamma_{g_i}) \ar[r] 
& \smash{\;\displaystyle\bigoplus_{j=1}^{\;n'}\;} \CM_k(\Gamma'_{g'_j}) \ar[r] 
& \smash{\;\displaystyle\bigoplus_{i=1}^n\;} \CM_k(\Gamma_{g_i}) \rlap{,} \\
(f_i)_{i=1}^n \ar@{|->}[r] & 
(f_{i_j}|_k\gamma_j)_{j=1}^{n'} \ar@{|->}[r] & 
\smash{\;\Bigl(\displaystyle\sum_{j,\gamma}f_{i_j}|_k\gamma_j\,|_k\gamma \Bigr){\strut}_{i=1}^n} \rlap{,} \\}$$
where, for each index~$i$, the sum extends over all pairs of indices $1\le j\le n'$ and elements $\gamma\in \GL_r(F)\cap g'_jKg_i^{-1}$ up to left multiplication by $\GL_r(F)\cap g'_jK'g_j^{\prime-1}$. Using the fact that $f_{i_j}|_k\gamma_j\,|_k\gamma = f_{i_j}|_k\,\gamma_j\gamma$ we can rewrite this as 
\UseTheoremCounterForNextEquation
\begin{equation}\label{HeckeDescriptionA}
(f_i)_{i=1}^n \ \longmapsto \Bigl(\sum_{j,\delta}f_{i_j}|_k\,\delta \Bigr){\strut}_{i=1}^n\strut_{\textstyle,}
\end{equation}
where, for each index~$i$, the sum extends over all pairs of indices $1\le j\le n'$ and elements $\delta\in \GL_r(F)\cap \gamma_jg'_jKg_i^{-1}$ up to left multiplication by $\GL_r(F)\cap \gamma_jg'_jK'g_j^{\prime-1}\gamma_j^{-1}$.

To analyse this sum note first that by construction we have $\gamma_jg'_j = g_{i_j}k_j^{-1}h$. For each $j$ the element $\delta$ therefore runs through $\GL_r(F)\cap g_{i_j}k_j^{-1}hKg_i^{-1}$ up to left multiplication by $\GL_r(F)\cap g_{i_j}k_j^{-1}hK'h^{-1}k_jg_{i_j}^{-1}$. 

For any $j$ and $\delta$ that occur in the sum this shows that $\delta g_i \in g_{i_j}k_j^{-1}hK$. Taking determinants and using the fact that $k_j\in K$ we deduce that $\det(g_i)$ and $\det(g_{i_j}h)$ represent the same coset in ${F^\times\backslash(\BAfF)^\times/\det(K)}$. The coset of $\det(g_{i_j})$ therefore depends only on $i$ and~$h$, but not on~$j$. By Proposition \ref{StrongApprox} it follows that $i_j$ depends only on $i$ and~$h$, but not on~$j$. For the rest of the proof we therefore fix indices $i$ and $i'$ such that $\det(g_i)$ and $\det(g_{i'}h)$ represent the same coset in ${F^\times\backslash(\BAfF)^\times/\det(K)}$, and we can restrict ourselves to indices $j$ with $i_j=i'$. 

Note that this already proves the last statement of the proposition. It also shows that $\delta$ lies in $\GL_r(F)\cap g_{i'}KhKg_i^{-1}$. Moreover, since $hK'h^{-1}<K$ and $k_j\in K$, we have $\GL_r(F)\cap g_{i'}k_j^{-1}hK'h^{-1}k_jg_{i'}^{-1} \subset \GL_r(F)\cap g_{i'}Kg_{i'}^{-1} = \Gamma_{g_{i'}}$. Thus any equivalence class of pairs $(j,\delta)$ in the sum (\ref{HeckeDescriptionA}) determines a unique coset $\Gamma_{g_{i'}}\delta$.

Suppose that two pairs $(j,\delta)$ and $(j',\delta')$ determine the same coset $\Gamma_{g_{i'}}\delta = \Gamma_{g_{i'}}\delta'$. Write $\delta'=\epsilon\delta$ with $\epsilon \in \Gamma_{g_{i'}}$. Since $\delta\in g_{i'}k_j^{-1}hKg_i^{-1}$ and $\delta'\in g_{i'}k_{j'}^{-1}hKg_i^{-1}$, it follows that $\delta'$ lies in both $\epsilon g_{i'}k_j^{-1}hKg_i^{-1}$ and $g_{i'}k_{j'}^{-1}hKg_i^{-1}$.
Multiplying by $g_i$ from the right we deduce that $\epsilon g_{i'}k_j^{-1}hk = g_{i'}k_{j'}^{-1}h$ for some $k\in K$. By the definition of $\Gamma_{g_{i'}}$ we have $g_{i'}^{-1}\epsilon^{-1}g_{i'}\in K$, and since $k_j$, $k_{j'}\in K$, we find that $k = h^{-1}k_jg_{i'}^{-1}\epsilon^{-1}g_{i'}k_{j'}^{-1}h
\in K\cap h^{-1}Kh = K'$. The calculation $\epsilon \gamma_jg'_jk = \epsilon g_{i'}k_j^{-1}hk = g_{i'}k_{j'}^{-1}h = \gamma_{j'}g'_{j'}$ now implies that $g'_j$ and $g'_{j'}$ represent the same double coset in $\GL_r(F)\backslash\GL_r(\BAfF)/K'$. By the choice of $g'_1,\ldots,g'_{n'}$ as representatives of these double cosets it follows that $j=j'$. Thus both $\delta$ and $\delta'$ lie in $\GL_r(F)\cap g_{i'}k_j^{-1}hKg_i^{-1}$, and hence $\epsilon = \delta'\delta^{-1}$ lies in $\GL_r(F)\cap g_{i'}k_j^{-1}hKh^{-1}k_jg_{i'}^{-1}$. 
Since also $\epsilon \in \Gamma_{g_{i'}} = \GL_r(F)\cap g_{i'}Kg_{i'}^{-1}$ and $k_j\in K$ and $hKh^{-1} \cap\nobreak K = hK'h^{-1}$, we then actually have $\epsilon\in \GL_r(F)\cap g_{i'}k_j^{-1}hK'h^{-1}k_jg_{i'}^{-1}$. 
This shows that the map sending an equivalence class of pairs $(j,\delta)$ in the sum (\ref{HeckeDescriptionA}) to the coset $\Gamma_{g_{i'}}\delta$ is injective.

Consider now an arbitrary element $\delta\in\GL_r(F)\cap g_{i'}KhKg_i^{-1}$. Choose $k\in K$ such that $\delta\in g_{i'}k^{-1}hKg_i^{-1}$. By the choice of $g'_1,\ldots,g'_{n'}$ there exists an index $j$ with $\GL_r(F)g_{i'}k^{-1}hK' = \GL_r(F)g_j'K'$. Since $\gamma_jg'_j = g_{i_j}k_j^{-1}h$, we deduce that $\GL_r(F)g_{i'}k^{-1}hK' = \GL_r(F)g_{i_j}k_j^{-1}hK'$. By the same argument as above it follows that $i'=i_j$, and we can find an element $\epsilon\in \GL_r(F)$ such that $\epsilon g_{i'}k^{-1}h \in g_{i'}k_j^{-1}hK'$. Since $hK'h^{-1}<K$ and $k_j$, $k\in K$, we then have $\epsilon \in \GL_r(F) \cap g_{i'}k_j^{-1}hK'h^{-1}k g_{i'}^{-1} < \GL_r(F)\cap g_{i'}Kg_{i'}^{-1} = \Gamma_{g_{i'}}$. Thus $\epsilon\delta \in \GL_r(F)\cap \epsilon g_{i'}k^{-1}hKg_i^{-1} = \GL_r(F)\cap g_{i'}k_j^{-1}hKg_i^{-1}$, and so the coset $\Gamma_{g_{i'}}\delta$ arises from the pair $(j,\epsilon\delta)$ in the sum (\ref{HeckeDescriptionA}). In other words the map sending an equivalence class of pairs $(j,\delta)$ in the sum (\ref{HeckeDescriptionA}) to the coset $\Gamma_{g_{i'}}\delta$ is surjective.

All this together shows that in (\ref{HeckeDescriptionA}) we can equivalently sum over all $\delta\in \GL_r(F)\cap g_{i'}KhKg_i^{-1}$ up to left multiplication by~$\Gamma_{g_{i'}}$. Also, since $f_{i_j} = f_{i'} \in \CM_k(\Gamma_{g_{i'}})$, the function $f_{i'}|_k\,\delta$ depends only on the coset $\Gamma_{g_{i'}}\delta$. This finishes the proof.
\end{Proof}

\medskip

Finally, we define Hecke operators on analytic Drinfeld modular forms as follows:

\begin{Def}\label{AnalyticHecke}
For any arithmetic subgroups $\Gamma$, $\Gamma'<\GL_r(F)$ and any element $\delta\in\GL_r(F)$ we define the associated \emph{Hecke operator} by
$$T_\delta:\ \CM_k(\Gamma') \longto \CM_k(\Gamma),\ 
f\longmapsto \sum\nolimits_{\gamma} f|_k\,\gamma,$$
where $\gamma$ runs through a set of representatives of the quotient $\Gamma'\backslash\Gamma'\delta\Gamma$.
\end{Def}

Using (\ref{fCocycle}) and Proposition \ref{Prop:Conj} one finds that this is well-defined, and by construction it depends only on the double coset $\Gamma'\delta\Gamma$. Also, since the action of $\GL_r(F)$ preserves cusp forms and $\CM_k(\Gamma) \cap \CS_k(\Gamma\cap\delta^{-1}\Gamma'\delta) = \CS_k(\Gamma)$, the Hecke operator induces a map
\UseTheoremCounterForNextEquation
\begin{equation}\label{HeckePreservesCuspForms}
T_\delta:\ \CS_k(\Gamma') \longto \CS_k(\Gamma).
\end{equation}

We can now rewrite the formula in Proposition \ref{HeckeDescription}
as follows.

\begin{Thm}\label{HeckeDescription2}
The map on the bottom in Proposition \ref{HeckeDescription} is equal to 
$$\qquad\xymatrix@C+20pt@R-15pt{
(f_i)_{i=1}^{n} \ar@{|->}[r] & 
\smash{\;\Bigl(\displaystyle\sum_{i',\delta} T_\delta(f_{i'}) \Bigr){\strut}_{i=1}^n} \rlap{,} \\}$$
where, for each index~$i$, the sum extends over all pairs of indices $1\le i'\le n$ and double cosets $\Gamma_{g_{i'}}\delta\Gamma_{g_i} \subset \GL_r(F)\cap g_{i'}KhKg_i^{-1}$. Again the index $i'$ that actually occurs depends only on $i$ and~$h$.
\end{Thm}

\begin{Proof}
By construction the set $\GL_r(F)\cap g_{i'}KhKg_i^{-1}$ is invariant under left multiplication by $\Gamma_{g_{i'}} = \GL_r(F)\cap g_{i'}Kg_{i'}^{-1}$ and right multiplication by $\Gamma_{g_i} = \GL_r(F)\cap g_iKg_i^{-1}$, and it is a finite disjoint union of double cosets $\Gamma_{g_{i'}}\delta\Gamma_{g_i}$. The formula results by direct computation from (\ref{fCocycle}).
\end{Proof}

\begin{Rem}\label{HeckeDescription3}
\rm In Theorem \ref{HeckeDescription2} it can happen that $\GL_r(F)\cap g_{i'}KhKg_i^{-1}$ decomposes into several double cosets. This is related to the fact that the algebraic Hecke operator $T_h$ is by construction defined over~$F$, whereas the analytic Hecke operator $T_\delta$ is only defined over~$\Cinf$. Thus if $M^r_{A,K \cap\nobreak h^{-1}Kh}(\Cinf)$ has more connected components than 
$M^r_{A,K}(\Cinf)$, their common field of definition $F_{K \cap\nobreak h^{-1}Kh}$ is a proper extension of the field of definition $F_K$ of the connected components of $M^r_{A,K}(\Cinf)$, and the algebraic Hecke operator $T_h$ can be viewed as an analytic Hecke operator $T_\delta$ followed by a trace map with respect to $F_{K \cap\nobreak h^{-1}Kh}/F_K$.
\end{Rem}

\begin{center}
\rule{8cm}{0.01cm}
\end{center}

\begin{minipage}[t]{5cm}{\small
Department of Mathematical Sciences \\
University of Stellenbosch \\
Stellenbosch, 7600 \\
South Africa \\
djbasson@sun.ac.za
}
\end{minipage}\hfill
\begin{minipage}[t]{5cm}{\small
School of Mathematical and Physical Sciences \\
University of Newcastle \\
Callaghan, 2308 \\
Australia \\
florian.breuer@newcastle.edu.au\\
 {\em and} \\
Department of Mathematical Sciences \\
University of Stellenbosch \\
Stellenbosch, 7600 \\
South Africa \\}
\end{minipage}\hfill
%
%
\begin{minipage}[t]{5cm}{\small
Department of Mathematics \\
ETH Z\"urich\\
8092 Z\"urich\\
Switzerland \\
pink@math.ethz.ch}
\end{minipage}

\end{document}